\documentclass[a4paper, 11pt]{article}
\usepackage{graphicx} 

\usepackage{comment}

\usepackage{todonotes}
\usepackage{enumitem}

\usepackage{xcolor}

\usepackage{amsmath}
\usepackage{amssymb}
\usepackage{amsthm}
\usepackage{bbm}

\usepackage{parskip}

\usepackage{tikz}

\usepackage{tikz-cd}
\theoremstyle{plain}
\newtheorem{theorem}{Theorem}

\newtheorem{proposition}[theorem]{Proposition}
\newtheorem{lemma}[theorem]{Lemma}
\newtheorem{corollary}[theorem]{Corollary}

\theoremstyle{remark}
\newtheorem{remark}[theorem]{Remark}

\theoremstyle{definition}
\newtheorem{definition}[theorem]{Definition}
\newtheorem{assumption}[theorem]{Assumption}

\numberwithin{equation}{section}
\numberwithin{theorem}{section}

\def\R{\mathbb{R}}

\def\C{\mathcal{C}}

\def\LL{\mathcal{L}}

\DeclareMathOperator{\ima}{Im}
\DeclareMathOperator{\tr}{tr}
\DeclareMathOperator{\sgn}{sgn}
\DeclareMathOperator{\rank}{rank}

\usepackage{mathtools}

\DeclareMathOperator\Div{div}
\DeclareMathOperator\Var{Var}

\usepackage{authblk}

\usepackage[T1]{fontenc}
\usepackage[utf8]{inputenc}
\usepackage[]{hyperref}
\usepackage[style=numeric, maxbibnames = 5]{biblatex}
\usepackage[a4paper, left=2.5cm, right=2.5cm, top=3cm, bottom=3cm]{geometry}
\title{A stochastic chemical system in contact with reservoirs}
\author[1]{E. Franco}
\author[1]{Jann Hendrik Schröders\footnote{\href{mailto:schroeders@iam.uni-bonn.de}{schroeders@iam.uni-bonn.de}, orcid ID: 0009-0006-6499-2653}}
\author[1]{J. J. L. Velázquez}
\affil[1]{Institute for Applied Mathematics, University Bonn}

\date{}

\begin{document}

\maketitle
\begin{abstract}
In this paper we analyse the behaviour of a non-reversible stochastic chemical system obtained as an approximation of a stochastic chemical system in which the concentration of some chemicals is higher than the concentration of the other chemicals. 
The master equation associated with this approximated system is characterised by a discrete part, modelling the reactions among chemicals that have low concentration, by a transport term and a diffusion term. The transport-diffusion terms describe the fluctuations in the concentration of abundant chemicals. The model contains a second order term is of order $\epsilon$, that describes the fluctuation in the concentrations of the abundant substances. We prove that the master equation has a spectral gap that is independent on $\epsilon$ and we analyse in detail the asymptotics of the solution to the master equation as $\epsilon \to 0$.
\end{abstract}

\section{Introduction}

 Biological systems are open systems that exchange free energy and chemicals with the environment in order to achieve their biological functions. We can imagine different ways in which the system exchanges chemicals and free energy with the environment. For example, chemicals can be injected and removed from the system via sources and sinks. Alternatively, the system can be assumed to be in contact with a reservoir of chemicals (chemostat). 
 A common way to model chemical systems in contact with reservoirs is to assume that the concentration of some chemicals in the system is constant in time. Examples of  chemicals that are often chemostated, i.e. assumed to have frozen concentrations, are energetic molecules, as ATP/ADP, GTP/GDP, NADH/NAD, glucose or phosphate groups. 
 
 The concentration of chemicals that are at a chemostated state is usually much higher than the one of non-chemostated chemicals that are present in the system. 
 This is the case, for example, in the kinetic proofreading mechanisms involved in immune recognition processes. Indeed in these systems few ligands interact with ATP, ADP and with receptors in order to trigger an immune response (see \cite{franco2026critical}). 
Also molecular motors, as for example kinesins, interact with high concentrations of ATP and ADP in order to move.
A similar situation occurs in DNA replication processes and protein synthesis, as well as many other biological processes (see for instance \cite{milo2015cell}).
 
If we have few chemicals of a given type in the system, it makes sense to model the number $X(t)$ of these chemicals as a stochastic process, namely a continuous-time Markov chain. 
Instead, the evolution in time of the concentration $N(t)$ of abundant chemicals can be modelled by deterministic systems, i.e. by systems of ODEs. These deterministic systems can be obtained starting form stochastic systems taking suitable scaling limits (see for instance \cite{anderson2015stochastic}). 
It is therefore natural to model chemical systems in which the chemostated chemicals are abundant while the other chemicals have low concentration using a hybrid stochastic-deterministic system (see \cite{ball2006asymptotic,remlein2025chemostat}).

In this paper we analyse in detail a master equation that approximates a chemical system in which a molecule interacts with abundant chemicals. The chemical reactions involving only the abundant chemicals are approximated by transport-diffusion terms in the master equation. The state of the molecule, instead, is described by a continuous-time Markov chain with a discrete state space. The concentrations of the abundant species is encoded in a continuous variable, that affects the transition rates of the Markov chain.  
The approximated master equation that we study is a model of a stochastic chemical system that is in contact with a reservoirs of chemicals, given by the abundant chemicals.

The model that we consider in this paper does not satisfy the detailed balance property.
The property of detailed balance is satisfied in a chemical system if every chemical reaction is balanced by the reverse reaction at steady state. It is a fundamental property that holds in every chemical system that is at constant temperature and that does not exchange chemicals and energy with the environment (see for instance \cite{polettini2014irreversible}).
Therefore chemical systems in which detailed balance fails are systems that are exchanging chemicals and free energy with the environment, for instance because they are in contact with a reservoir of chemicals. 
In \cite{franco2025detailed} we prove that if we freeze at constant values the concentrations of some chemicals in a deterministic system in which detailed balance holds, we obtain a chemical system, involving only the non-frozen chemicals, that does not satisfy the detailed balance property, unless the reaction rates are fine tuned, or the frozen concentrations are chosen at equilibrium values, or the chemical system satisfy a very particular topological condition. 

In \cite{franco2026reservoir,franco2026flux} we adapt the results obtained in \cite{franco2025detailed} to the stochastic setting. 
More precisely, in \cite{franco2026flux} we prove that adding sources and sinks to a stochastic chemical system in which detailed balance holds, we obtain a system in which detailed balance fails, unless a topological condition is satisfied or the chemical rates are fine tuned. 
In \cite{franco2026reservoir} we derive, using formal arguments, a stochastic chemical system in which detailed balance fails starting from a stochastic chemical system in which detailed balance holds and in which some chemicals have high concentration. 
The chemicals that have high concentrations act as a reservoir for the stochastic system that describes the dynamics of other chemicals.
Since detailed balance holds for the larger stochastic chemical system that includes the abundant chemicals we expect the system to converge to the Gibbs free distribution as time tends to infinity. As a consequence this reservoir keeps the microscopic system out of equilibrium only for times of order one, i.e. it is a transient reservoir. 

In \cite{franco2026reservoir} we also derive a class of stochastic chemical systems in which detailed balance fails starting from a stochastic chemical system in which sources and sinks keep the concentrations of some chemicals high. 
The rigorous analysis of the master equation associated with this non-reversible stochastic chemical system is challenging and is not discussed in  \cite{franco2026reservoir}. 
The model that we study in this paper is a system in which detailed balance fails and that we derive starting from a chemical system in which detailed balance fails and the concentrations of some chemicals are high. The initial system, that includes the abundant chemicals, could be obtained adding sources and sinks to a system in which detailed balance holds. 
However in this paper we do not do this, in order to keep the model simple enough to be able to obtain clean mathematical results.

 The model that we consider can be used to describe several chemical systems that work in out of equilibrium conditions. For example it can be used to model of kinesines, DNA replication and protein synthesis. By definiteness in this introduction we will interpret it as an idealised model of a molecular motor/kinesin.
Motor proteins, such as the different myosins and kinesins, in living cells perform crucial tasks, from transporting organelles and vesicles through the cell to muscle contraction.
For example, kinesin molecules can perform a ``walking'' movement along a filament by hydrolysing one ATP per ``step'' it performs, \cite{kinesin_1997}, where the binding of ATP to the kinesin gives the stepping process its direction \cite{structural_1999}. This energy conversion mechanism from ATP to protein motion, which is shared among kinesin and myosin, has allowed for the evolution of a vast number of motor proteins of different functions, see e.g., the reviews \cite{way_2000} and \cite{kinesin_2009}.

The class of  models that we consider is similar to the ones described in \cite{walking-molecules}. 
We consider a chemical system consisting of a kinesin $K$ that can be at $n $ different sites and is surrounded by solution of energetic molecules, as for instance ATP and ADP molecules. 
The energetic molecules are the abundant chemicals in the system and are used by the molecule $K$ to change state/position.
Even though we will analyse a more general chemical system, the types of reactions that we have in mind are chemical reactions of the following form  
\begin{equation}\label{network}
\begin{tikzcd}
	{(ATP)} & {(ADP)} && {(ATP)+(k)} & {(ADP)+(k+1)}
	\arrow["{\rho_0}", shift left=2, from=1-1, to=1-2]
	\arrow["{-\rho_0}", shift left=2, from=1-2, to=1-1]
	\arrow["{\rho_k}", shift left=2, from=1-4, to=1-5]
	\arrow["{-\rho_k}", shift left=2, from=1-5, to=1-4]
\end{tikzcd}
\end{equation} 
where $(k)$ is the state of the molecule $K$. 
Since we are interested in studying active systems in which detailed balance fails we assume that $1 =1+n$. In this way we introduce a cycle in the system. 

Typically in biological applications molecular walkers move along straight filaments. When they reach the end of the filaments they fall, enter a diffusive state and then go back to the filament. We could think about the state $n-1 $ as the state in which the motor falls and as $n $ as the diffusive state. After the diffusive state $n$ the motor goes back to the filament. 
Notice that there are also artificial motors, as for instance catenated wheel motors, that sometimes are artificially forced to walk on rings (see \cite{leigh2003unidirectional}).

 We associate to this chemical system the master equation whose solution is the probability $F(t, n_{ATP}, n_{ADP}, k ) $ that at time $t>0 $ the system is at state $(n_{ATP}, n_{ADP}, k ) $ where $n_{ATP} $ and $n_{ADP}$ are the number of $ATP$ and of type $ATP$ while $k $ is the state of $K  $. 
As we illustrate in Section \ref{Sec_Motivation}, when the concentration of chemicals of type $ATP$ and $ADP$ is high and the chemical rates satisfy suitable scaling limits, it is possible to approximate  the master equation associated with \eqref{network}, for times of order $1$, with the following system of equations
\begin{equation}\label{intro_disc_cont_eq}
   \begin{aligned}
    \partial_t f_k (t, x)  &=  \epsilon\partial_x^2 f_k (t, x)   +\partial_x (x f_k (t, x)  ) +(a_{k-1} (x)  f_{k-1} (t, x)  -a_k (x) f_k (t, x)  ) \\
    & +(b_{k+1} (x)   f_{k+1} (t, x)  - b_k(x)   f_k (t, x)  ),\quad  k=1,\dots,n ,  \quad x \in \mathbb R. 
\end{aligned} 
\end{equation}

Here $f_k (t, x ) $ is the probability that $K$  is at state $ k $ at time $ t >0 $ while the concentrations of $ATP$ and $ADP$ are close to equilibrium. 
The variable $x\in\R$ encodes the fluctuation in the concentrations of $ATP$ and $ADP$ away from the equilibrium. 
The parameters $a_k (x) $ and $b_k (x) $ are the chemical rates associated with the chemical reactions $\rho_k  $ and $ \rho_{ - k } $ and they depend on the fluctuations in the concentrations of $ATP$ and $ADP$.

The parameter $\epsilon \ll  1 $ in equation \eqref{intro_disc_cont_eq} accounts for the fact that we are approximating the master equation of a stochastic process with discrete state space with a continuous model. This is justified by the assumption $n_{ ATP} , n_{ADP} \gg 1 $ and the assumption that the rates of the reactions $\rho_k, - \rho_{  k }  $ are much slower than the rates of the reactions $\rho_0,- \rho_0 $. In particular the terms $ \epsilon\partial_x^2 f_k (t, x)   +\partial_x (x f_k (t, x)  )$ approximate the chemical reactions $\rho_0$ and $ - \rho_0 $. 

Notice that the approximation \eqref{intro_disc_cont_eq} is valid only for times that are of order one and that eventually we expect the chemical system \eqref{network} to converge to a stationary distribution. 
In \cite{franco2026reservoir} we study in detail how to derive the analogous to the approximating equation \eqref{intro_disc_cont_eq} for a class of models that is more general than the one in \eqref{network}.
Moreover in that paper we also analyse in detail the behaviour for large times, or more precisely for times of order $\epsilon^{-1 }$. 
 
An important feature of the effective model \eqref{intro_disc_cont_eq}, that we study in this paper, is the fact that the chemical rates $a_k(x)  , b_k (x) $ depend on the concentration of $ATP$ and $ADP. $ This dependence can be observed also in biological systems. 
In \cite{kinesin_1997} the speed of the kinesin walk increases with the supply of ATP. 
In addition, artificial molecular walkers, whose movements can be controlled by changing the chemical composition of the solution in which they are contained, have been constructed (see for example \cite{clocked_2026}). 

In this paper we analyse in detail a class of equations that include  \eqref{intro_disc_cont_eq} as a special case, i.e. we study  the system
\begin{equation}\label{intro: gen disc_cont_eq}
   \begin{aligned}
    \partial_t f_k(t, x)  &=\epsilon\partial_x^2 f_k(t, x)  +\partial_x(H'(x)  f_k (t, x) ) + \mathbf{e}_k^\dagger  M(x) f (t, x) ,\quad k=1,\dots,n 
\end{aligned} 
\end{equation}
where $H $ is a suitable convex function, $M(x) $ is a matrix that describes a chemical system and $f=\sum_{ k =1}^n f_k \mathbf{e}_k $. 
We first of all prove the existence of a unique stationary solution to equation \eqref{intro: gen disc_cont_eq}. To this end we analyse the log-time behaviour of the dual to equation \eqref{intro: gen disc_cont_eq}. 
In order to do this we use the theory of Lyapunov equations (see \cite{stein1952lyapunov}) as well as a Poincaré inequality due to \cite{bakryemry}. This allows to prove the existence of a spectral gap for the dual equation. The spectral gap turns out to be independent on the parameter $\epsilon $. 
This allows us to prove that for fixed $\epsilon >0 $, the time dependent solution $ f $ converges to the unique stationary solution $\bar f_\epsilon$ as $t \to \infty $, with an exponential rate that does not dependent on $\epsilon$, i.e. 
 \begin{equation} \label{intro:convergence}
     \vert\vert f(t)-\bar{f}_\epsilon \vert\vert_{L^2(d\mu)}\lesssim \exp(-\lambda t) \quad t >0 
 \end{equation}
where $\mu (dx) = \gamma (x )  dx $ with $ Z_\epsilon\eqqcolon \left( \int_\R \exp(-H/\epsilon)\,dx \right)^{-1} $ and $\gamma (x) \eqqcolon \frac{1}{Z_\epsilon}\exp(-H(x) /\epsilon )$.

 We consider convex potentials $H$ such that $H'(0)=0$, as a consequence there is a boundary layer at $x =0$. 
Therefore, in order to study the behaviour of the stationary solution as $\epsilon \to 0 $ we consider the Ansatz 
\begin{equation*}
    f_k (t, x ) = \rho_k (t, x) e^{\frac{ - H (x)}{\epsilon}}  \quad k \in \mathbb N, \quad x \in \mathbb R, \quad t >0. 
\end{equation*}
and study the behaviour of $\rho_\epsilon := \sum_{ k =1}^n \rho_k\,\mathbf{e}_k $ as $ t \to \infty $ and as $\epsilon \to \infty$.
It is easy to see that the convergence \eqref{intro:convergence} implies that $\rho_\epsilon$ converges to a stationary profile, i.e. for fixed $\epsilon >0$, $\rho (t) \rightarrow g_\epsilon $ as $ t \to \infty $ where $g_\epsilon $ is the unique solution to 
\begin{equation}\label{intro:g-eq}
    0 =\epsilon\frac{d^2}{dx} g_{\epsilon,k }-H'(x)\,\frac{d}{dx} g_{\epsilon,k} + e_k^\dagger M(x)g_\epsilon.
\end{equation}

We then analyse in detail the limit of $g_\epsilon $ as $\epsilon \to 0 $. 
We deduce that $ g_\epsilon \to g $ globally uniformly on $\R $ where $g$ is the solution to the limit equation
\begin{equation} \label{intro_asympt_eq} 
    g'(x) H' (x)  =  M(x) g(x). 
\end{equation}
In order to prove this we analyse in detail the dynamics around the boundary layer $x=0$. 

Notice that since the model that we consider does not satisfy, in general, the detailed balance property, therefore the stationary solution to equation \eqref{intro: gen disc_cont_eq} is a non-equilibrium stationary solution.
To see this, it is convenient to rewrite the stationary equation \eqref{intro:g-eq} in flux form, i.e. as 
\[ 
0=\Div(J( \bar f_\epsilon) )+ F^M( \bar f_\epsilon ) \mathbf{1}  
\]
where 
\begin{equation} \label{cont flux}
J(\bar f_\epsilon)=\epsilon \partial_x f + H'\,f =\epsilon\partial_x g_\epsilon \,\gamma
\end{equation} 
is the flux in the $x$-variable while the matrix $F^M$, defined in \eqref{matrix_flux}, encodes the net conversion rates between the $n$ states. After choosing an orientation for the reactions, we define the discrete flux, at the state $(x,k)$ as the vector
\begin{equation}\label{disc_flux}
    I(x,k) \coloneqq  \sum_{\ell\not= k} \sgn(\ell-k)(M_{k,\ell}(x)\bar{f}_k(x)-M_{\ell,k}(x)\bar f_\ell(x) )\mathbf{e}_\ell.
\end{equation}
Detailed balance holds if and only if both the discrete and the continuous fluxes are zero for the stationary states i.e. if $F_{ M(x)}\,(\bar f_\epsilon ) =0 $ and $J (\bar f_\epsilon ) =0$. 
We prove that detailed balance holds only if $ M(x)= M $ for every $x \in \mathbb R $. Moreover we prove that if detailed balance does not hold, then $J (\bar f_\epsilon ) \approx o_\epsilon(1)\,\delta_0 $.
In the limit as $\epsilon \to 0 $ we have that $J (\bar f_\epsilon ) \to 0  $ and $ \bar f_\epsilon \to \delta_0$.  Therefore the $F_{ M(x)}  (\bar f_\epsilon )$ is supported only in $x=0 $ in the limit as $\epsilon \to 0$. The detailed balance property for the limiting system holds only if the chemical system associated with $M(0)$ satisfies the detailed balance property. Since we assume that the original system \eqref{network} does not satisfy the detailed balance property, this means that the limiting function $g$ is a stationary solution with non-zero discrete fluxes.  

In order to get more intuition let us assume that the matrix $M$ in \eqref{intro: gen disc_cont_eq} is the one associated with the model of kinesin \eqref{intro_disc_cont_eq}.
In the limit as $\epsilon \to 0 $ the matrix of the fluxes have the following form 
\begin{equation} \label{constant fluxes}
  I(x,k)\approx \left(J \mathbf{e}_{k+1} - J\mathbf{e}_k\right)\delta_0.
\end{equation}
for a constant $J \neq 0$. In other words the stationary solution $g$ is a stationary solution with constant fluxes. This stationary solution with fluxes describes a molecule that walks on the filament, i.e. on the sites $\{ 1, \dots, n \}$. The sign of $J$ determines the direction of this motion.

Models that have mathematical analogies and also important differences with  \eqref{intro_disc_cont_eq} are the models of molecular motor analysed in \cite{perthame2009asymmetric} and in \cite{chipot2004transport}. These models are inspired by the model of Brownian ratchet in \cite{feynman1963feynman}. 
As \eqref{intro_disc_cont_eq} the models considered in \cite{chipot2004transport,perthame2009asymmetric} are Fokker-Planck systems coupled with discrete dynamics induced by linear chemical reactions.  
Despite the similarities between the equations analysed in \cite{chipot2004transport,perthame2009asymmetric} and the equation that we analyse in this paper, the model that we consider here differs from the models studied in  \cite{chipot2004transport,perthame2009asymmetric}.
Our model is derived starting from a stochastic chemical systems in which some chemicals have high concentrations. The transport-diffusion terms in our model describe the chemical reactions among abundant chemicals. 
Accordingly, the continuous variable describes the evolution of the concentration of the abundant chemicals. Instead the discrete states $\{1, \dots, n \} $ are the sites along the filaments and represent space. 
The continuous variable in \cite{chipot2004transport,perthame2009asymmetric} describes the space on which the motor walks. The diffusion-transport term describe the movement of the molecule. The discrete variables in \cite{chipot2004transport,perthame2009asymmetric} describe the chemical state of the motor.

One of the goals in \cite{chipot2004transport,perthame2009asymmetric} is to prove that the model induces a directionality in the walking of the motor. In order to achieve this directionality the authors consider periodic asymmetric potentials. The role of this particular choice of potentials is to break the detailed balance property. Indeed we cannot expect to have any directionality in systems in which detailed balance holds (see \cite{feynman1963feynman}).  
In our paper we consider simpler potentials, that are independent on the discrete state $(k)$, describing the position of the kinesin, and that are convex. We obtain the desired directionality from the assumption that detailed balance fails for our chemical system \eqref{network}. This guarantees the existence of constant flux solutions $g $ (see \ref{constant fluxes}). The sign of the constant fluxes determines the directionality.  We stress that these are constant flux solutions in the discrete variable, that for us is space. 

We conclude the introduction with a comment on our choice of potentials $H$. 
In this paper we consider mainly convex potentials $H$. As we will see in Section \ref{Sec_Motivation} this is the natural assumption for the model that we consider in this paper.  

\subsection{Plan of the paper}
The plan of the paper is the following. 
In Section \ref{sec_model} we introduce the model and present the main results. 
In Section \ref{Sec_Motivation} we derive the model that we study in this paper starting from a chemical system in which detailed balance fails and in which some chemicals are abundant. 
In Section \ref{sec_well_posed} we prove the existence of a solution to the time dependent master equation \eqref{intro: gen disc_cont_eq} and the existence of a unique solution to the stationary equation \eqref{intro:g-eq}. 
We then state some useful results on the matrix $M(x) $ in Section \ref{sec_M}. 
These results allows us to study in detail the behaviour of the solution to equation \eqref{intro_asympt_eq} in Section \ref{sec_asymptotics}. 
In Section \ref{sec_convergence} we prove the convergence of $g_\epsilon $ to $g$ as $\epsilon \to 0$. 
Finally in Section \ref{sec_fluxes} we compute the constant fluxes associated with $g_\epsilon $ and with its limit $g$.

\subsection{Notation} 
 We write $X^\dagger $ to denote the transpose of a matrix $X$ and $A^*$ to denote the (formal) adjoint of a linear operator $A$. $I_d\in \R^{n\times n}$ denotes the identity matrix. All vectors $\mathbf{v} \in \R^n$ are column vectors and $(\mathbf{e}_k)_{k=1}^n$ denotes the standard basis of $\R^n$. We use $\langle v,w\rangle$ to denote the inner product of two elements of a Hilbert space, but for finite dimensional vector spaces we also simply write $\mathbf{v}\cdot \mathbf{w}$. In chains of inequalities, a generic constant $C>0$ is allowed to change in every step.
We use the notation $f\lesssim g $ when there exists a constant $C>0$ such that $ f \leq C g $. 
We use the notation $f \approx g $ to indicate that $f$ and $g $ are roughly of the same order.

\section{Model and main results}\label{sec_model}

The class of models that we consider in this paper is the following. 

On $\R_+\times \R$, we consider the system of $n$ PDEs given by
\begin{equation}\label{disc_cont_eq}
   \begin{aligned}
    \partial_t f_k(t, x)  &=\epsilon\partial_x^2 f_k(t, x)  +\partial_x(H'(x)  f_k (t, x) ) + \mathbf{e}_k^\dagger  M(x) f (t, x) ,\quad k=1,\dots,n\\
    \partial_x f_k(t,0)&=0,\\
\end{aligned} 
\end{equation}
with the convention that $n+k\equiv k$ and where we denote as $f$ the vector of the probabilities $( f_k )_{k =1}^n $, i.e. 
\[ 
  f \coloneqq \sum_{k=1}^n f_k \mathbf{e}_k. 
\] 

We make the following assumptions on the function $H $ and on the matrix $M(x) $. 
First of all we need some  boundedness and regularity conditions on the matrix $M(x)$, which we will need later on to be able to talk about the asymptotic behaviour of the solution, especially near $\pm \infty$. 
\begin{figure}
\centering
 \begin{tikzcd}
	&& {} & \dots & {} & {} & \dots & {} & {} & \\
	& {} & \dots & {(i)} &&& {(j)} &&& \dots \\
	&& {} && {} & {} && {} & {} \\
	& {} & {} && {} & {} && {} \\
	\dots &&& {(\ell)} &&& {(k)} & \dots \\
	& {} & {} & \dots & {} & {} & \dots
	\arrow[shift left=2, from=1-4, to=2-4]
	\arrow["{\scriptsize{(A)}}", shift right=3, bend right=60, from=1-5, to=1-6]
	\arrow[shift left=2, from=1-7, to=2-7]
	\arrow["{\scriptsize{(A)}}", shift right=3, bend right=60, from=1-8, to=1-9]
	\arrow[shift left=2, from=2-3, to=2-4]
	\arrow[shift left=2, from=2-4, to=1-4]
	\arrow[shift left=2, from=2-4, to=2-3]
	\arrow[shift left=2, from=2-4, to=2-7]
	\arrow[shift left=2, from=2-4, to=5-4]
	\arrow[shift left=2, from=2-7, to=1-7]
	\arrow[shift left=2, from=2-7, to=2-4]
	\arrow[shift left=2, from=2-7, to=2-10]
	\arrow[shift left=2, from=2-7, to=5-7]
	\arrow[shift left=2, from=2-10, to=2-7]
	\arrow["{\scriptsize{(A)}}"{pos=0.1}, shift right=5, bend right=60, from=3-5, to=4-5]
	\arrow["{\scriptsize{(B)}}", shift right=5, bend right=60, from=3-6, to=3-5]
	\arrow["{\scriptsize{(A)}}"{pos=0.1}, shift right=5, bend right=60, from=3-8, to=4-8]
	\arrow["{\scriptsize{(B)}}", shift right=5, bend right=60, from=3-9, to=3-8]
	\arrow["{\scriptsize{(A)}}", shift right=5, bend right=60, from=4-2, to=4-3]
	\arrow["{\scriptsize{(B)}}"{pos=0.1}, shift right=5, bend right=60, from=4-3, to=3-3]
	\arrow["{\scriptsize{(A)}}", shift right=5, bend right=60, from=4-5, to=4-6]
	\arrow["{\scriptsize{(B)}}"{pos=0.1}, shift right=5, bend right=60, from=4-6, to=3-6]
	\arrow[shift left=2, from=5-1, to=5-4]
	\arrow[shift left=2, from=5-4, to=2-4]
	\arrow[shift left=2, from=5-4, to=5-1]
	\arrow[shift left=2, from=5-4, to=5-7]
	\arrow[shift left=2, from=5-4, to=6-4]
	\arrow[shift left=2, from=5-7, to=2-7]
	\arrow[shift left=2, from=5-7, to=5-4]
	\arrow[shift left=2, from=5-7, to=5-8]
	\arrow[shift left=2, from=5-7, to=6-7]
	\arrow[shift left=2, from=5-8, to=5-7]
	\arrow["{\scriptsize{(B)}}", shift right=5, bend right=60, from=6-3, to=6-2]
	\arrow[shift left=2, from=6-4, to=5-4]
	\arrow["{\scriptsize{(B)}}", shift right=5, bend right=60, from=6-6, to=6-5]
	\arrow[shift left=2, from=6-7, to=5-7]
\end{tikzcd}
\caption{\footnotesize{A chemical network of the type considered in Equation \eqref{disc_cont_eq}: The transitions between the different states take place by consumption of the macroscopic substances $A$ and $B$, which are supplied by the reservoirs }}\label{fig:chemical}
\end{figure}
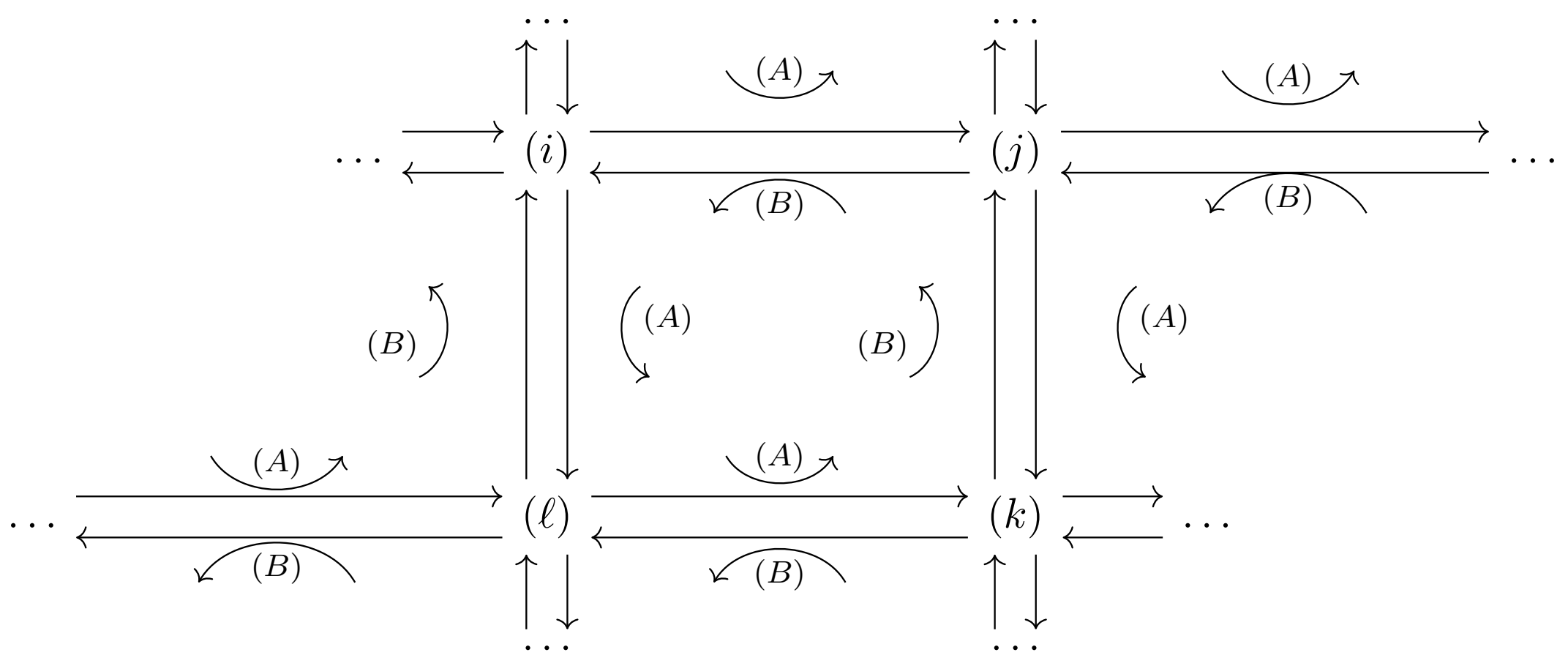
\begin{assumption}[Assumptions on the matrix $M(x)$]  \label{assumption_M_structure} \hfill

We assume that $M:\R\to \R^{n\times n}$ is smooth and component-wise bounded. Moreover we assume that 
\begin{enumerate}
    \item for every $ x \in \mathbb R$ we have that $M_{ij} (x) \geq 0 $ whenever $i\not=j$ and $\sum_{i=1}^n M_{ij}(x) =0$, i.e., $M_{ii}(x) =-\sum_{j=1}^n  M_{ij} (x) $, for all $j \in \{ 1, \dots N \}$. 
    \item\label{item_rank} $\rank M(x)=n-1$ for every $x\in \mathbb R$. 
    \item If $\lambda_k(x) \not =0$ are the non-zero eigenvalues of $M(x)$ then $\min_{k,x} \vert \Re \lambda_k(x)\vert \ge c_m>0$.
    \item We assume also that the limits
\begin{equation}\label{limits_M}
    \lim_{x\to \pm\infty} M(x)=M(\pm\infty)
\end{equation}
exist and that $r(x)\coloneqq M(x)-M(\pm\infty)=O(|x|^{-\alpha})$ as $|x| \to \infty $ with $\alpha>0$. 
\end{enumerate}
\end{assumption}
In view of Section \ref{Sec_Motivation}, where we will derive  the potential $H(x)=x^2/2$, we make the following assumption.
\begin{assumption}[Assumptions on the function $H$] \label{assumption_H}
We assume $H:\R\to\R$ to be smooth, uniformly convex, $0<c\le H''(x)\le C$, and bounded from below. 

Without loss of generality, we therefore assume that $H(0)$ is the unique minimum of $H$.
\end{assumption} 
\begin{remark}\label{rem:assumption_H}
Convexity is not necessary, and, in fact, all results will remain valid, if we make the following assumptions, which allow for non-convex potentials $H$. Assume $H:\R\to \R$ to be a smooth Morse function and to satisfy the following inequalities with $C,\lambda>0$.
\begin{align*}
    x^2/2 -C&\le H(x), \quad \forall x \in \R \\
    \sup_{x \in \mathbb R} \vert H'(x)- x\vert&\le \lambda,\\
    \sup_{ x\in \mathbb R}\vert H''(x)\vert&\le \lambda
\end{align*}
Further, we assume that $H$ is also bounded from above by some polynomial.
\end{remark}
We now make some remarks on Assumption \ref{assumption_M_structure}. 
\begin{remark}
For example, let us define $M :  \R\to \R^{n\times n} $ as follows
   \begin{align}
    M(x)\,v  &\coloneqq \sum_k [(a_{k-1}v_{k-1}-a_k(x) v_k)+(b_{k+1}(x) v_{k+1}- b_k(x) v_k)] \mathbf{e}_k, \quad  \forall v \in \R^n. 
\end{align}  
Moreover, assume that the $a_k,\,b_k$ are smooth, bounded and  uniformly positive, that is, there exists $c_M$, such that
\begin{equation*}
    \min_{k,x} \{a_k(x),b_k(x)\}\ge c_M>0.
\end{equation*} 
Then $M $ satisfies Assumption \ref{assumption_M_structure}.
\end{remark}
\begin{remark}\label{rem_net}
We interpret the matrix $M(x)$ as a matrix of linear chemical reactions, i.e. the columns of $M(x)$. More precisely, for some $i \neq j \in \{ 1, \dots, n \}   $ the chemical reaction \begin{equation*}
    \begin{tikzcd}
	{(i)} & {(j)}
	\arrow[shift left=0, from=1-1, to=1-2]
\end{tikzcd}
\end{equation*} takes place at rates $M_{i,j}(x)$, see Figure \ref{fig:chemical}. Therefore, the point \ref{item_rank} of Assumption \ref{assumption_M_structure} is satisfied if the graph associated with the chemical reactions in $M(x)$ is connected.
\end{remark}
Notice that point 1. in Assumption \ref{assumption_M_structure} implies that $(1,\dots,1)^\dagger  M(x)=0$. Therefore a formal computation suggests  that we should expect to have mass conservation, i.e. 
\begin{equation*}
    \frac{d}{dt} m(t)=\int_\R \sum_{k=1}^n \partial_t f_k(t, x) \,dx=0.
\end{equation*}
\par
Since system \eqref{disc_cont_eq} is a weakly coupled system of Fokker-Planck equations, it is natural to make the Ansatz  \[ 
f_k(t, x) =  c_\epsilon \rho_k(t, x)\exp(-H(x) /\epsilon), 
\]
that consists in writing $f_k$ in terms of a density $\rho_k $ multiplied against the Gibbs measure associated to  $H'$.
This leads to the following system for $\rho_k$.  
\begin{align}\label{disc_cont_eq rho}
      f_k (t, x) &\coloneqq \exp(-H(x)/\epsilon)\,\rho_k (t, x)\\
    \partial_t  \rho_k (t, x)  &=\epsilon\partial_x^2 \rho_k (t, x) -H'(x)\,\partial_x  \rho_k (t, x) +\mathbf{e}_k^\dagger  M(x) \rho (t, x),\, k=1,\dots,N \nonumber
\end{align}
where $\rho= \sum_{ k =1 }^n \rho_k \mathbf{e}_k $. 
Since, as explained above the system \eqref{disc_cont_eq} is such that $ \dfrac{ d }{dt } m(t)  =0 $ it is natural to check the integrability of the Gibbs measure.
To this end we notice that Assumption \ref{assumption_H}  implies that the function $\exp(-H/\epsilon)$ is integrable on $\R$. Therefore, we can define 
\begin{equation}
    \gamma(x)\coloneqq\left( \int_\R \exp(-H/\epsilon)\,dx \right)^{-1} \exp(-H(x)/\epsilon)\eqqcolon \frac{1}{Z_\epsilon}\exp(-H/\epsilon).
\end{equation}
We define $\mu(dx) = \gamma (x) dx $. 

\subsection{Main results}
In this paper we study the long-time behaviour of the solution to equation \eqref{disc_cont_eq} and its limit as $\epsilon \to 0$. First of all let us mention that the existence of a time dependent solution and of a stationary solution can be proven using classical techniques (see Proposition \ref{prop:existence time dep} and Lemma \ref{lem:existence stationary}). 
The first main of the paper is the existence of a spectral gap that does not depend on $\epsilon$ for equation \eqref{disc_cont_eq}. More precise we prove the following theorem. 
 \begin{theorem} \label{thm:main reults spectral gap}
 Let $\rho \in L^\infty_{loc}(L^2(d\mu))$ be a solution to \eqref{disc_cont_eq rho} with initial mass $m:= \int_\R \sum_{k=1}^n  \rho_k(0, x) \,dx $. Let $g_\epsilon$ be the unique stationary solution to \eqref{disc_cont_eq rho} with mass $m $. Then
 \begin{equation*}
     \vert\vert \rho(t)-g_\epsilon\vert\vert_{L^2(d\mu)}\le C\exp(-\lambda t), \quad t >0 
 \end{equation*}
 where $\lambda >0$ does not depend on $\epsilon$ and $C>0$. 
 \end{theorem}  
 As anticipated above in order to prove this theorem we use a Poincarè estimate as well as the theory of Lyapunov equations.

 We then analyse the asymptotics of $g_\epsilon $ as $\epsilon \to 0 $. We prove that, as expected, it converges to the equation that we obtain taking formally the limit as $\epsilon \to 0 $ in \eqref{intro:g-eq}.  
\begin{theorem} \label{thm: convergence as epsilon to zero} 
 As $\epsilon\to0$, the solutions $g_\epsilon=(g_{\epsilon,k})_{k=1}^n$ to equation \eqref{intro:g-eq} converge \emph{globally uniformly on $\R$} to the unique positive solution $g=(g_k)_{k=1}^n \in C_b^1(\mathbb R) $ satisfying
 \begin{equation*}
      \begin{cases}
     H'(x) \frac{d}{dx} g(x)&=M(x)g(x),\\
     \sum_{k=1}^n g_{\epsilon,k} &=\sum_{k=1}^n  g_k.
 \end{cases}
 \end{equation*}
 \end{theorem} 
 In order to prove this result we proceed with two steps. First of all we prove that up to a subsequence 
 \[ 
 g_\epsilon \to g \  \text{ locally uniformly in } \  \mathbb R \setminus \{ 0 \} \text{  as }  \epsilon \to 0.
 \] 
 Then we analyse in detail the behaviour of $g_\epsilon $ in $ x =0 $. 
 To this end we find explicit asymptotics for $g_\epsilon $ at $x=0$. 
To do this we rewrite system \eqref{intro:g-eq} as a system of equations of the form 
 \[
 \epsilon g_{\epsilon,k}'' - x g_{\epsilon, k }  + \lambda_k g_{\epsilon, k } =0. 
 \] 
 These equations can be written as Weber's equation via a change of variables. This allows us to study in detail the asymptotics of $g_\epsilon $ as $ x \to \pm \infty $ and around $x=0$ using parabolic cylinder functions. 
Using these explicit solutions we deduce that for $ |x | \leq \delta $, with $\delta $ sufficiently small, the following approximation holds 
\[
\| g_\epsilon (x)  -\bar  \ell \| =O(\epsilon^{1/4}) \quad  \text{ as } \epsilon \to 0. 
\] 
 where $ \bar \ell $ is the unique vector such that $ \bar \ell  \in \ker M(0)$ and $m =\sum_{k=1}^n \ell_k  $. Analysing the system $H' g' = M g $ (see Proposition \ref{prop: existence and uniqueness solution limiting problem}) we deduce that, as expected, we must have $g(0)= \bar \ell $.
 Combining this result with the uniform convergence  $g_\epsilon \to g $ we prove Theorem \eqref{thm: convergence as epsilon to zero}. 
 
One of the main features of systems in which detailed balance fail is the existence of constant flux stationary solutions.
For these systems we introduce the following matrix that encodes the net conversion rate from $i$ to $j$
\begin{equation}\label{matrix_flux}
    (F^M\,v)_{i,j}\coloneqq \begin{cases}
       m_{i,j}v_i-m_{j,i}v_j,\quad &i\not=j\\
       0,\quad &i=j.
    \end{cases} 
\end{equation}
We study in detail the fluxes associated with the stationary solution $ \bar f_\epsilon $ and prove the following statement. 
\begin{theorem} \label{thm:fluxes_main results}
Fix $\epsilon>0$. Then the continuous flux $J(\overline f_\epsilon )$ given by \eqref{cont flux} is zero if and only if $M(x)\equiv \tilde{M}\in \R^{n\times n}$ and $\ker \tilde{M}=\ker F_{\tilde{M}}\not=\{0\}$ where $F_{\tilde M }  $ is the matrix  associated with the fluxes, defined in \ref{matrix_flux}. 
Moreover the continuous flux $J(\bar f_\epsilon ) $ vanishes as $\epsilon\to0$ as a distribution and 
\begin{equation}
 F^{M(x) } (\bar f_\epsilon ) =  \delta_0 (x)  F_{M(0) } (g(0) ) \text{ as } \epsilon \to 0
\end{equation}
and
\begin{equation*}
    \mathcal{J}(x,k)\approx -\delta_0\sum_{k\not=\ell} \sgn(\ell-k) F^{M(0)}_{k,\ell}(g(0))(0,\mathbf{e}_k)^\dagger,\quad \text{as } \epsilon\to 0.
\end{equation*}
\end{theorem}
 Notice that for fixed $\epsilon >0$ the statement above guarantees that if $M(x)$ is non-constant, the reservoirs of $A$ and $B$ keep the system \eqref{intro: gen disc_cont_eq} out of equilibrium. 
 In the limit as $\epsilon \to 0 $ the continuous flux vanishes and the discrete fluxes are zero for stationary solutions if and only if detailed balance holds for system \eqref{network}. 
 
\section{Derivation of the model}\label{Sec_Motivation}
The goal of this section is to explain why we are interested in analysing equation \eqref{disc_cont_eq}. To this end we will show that it is possible to derive equation \eqref{disc_cont_eq} starting from the master equation that models a stochastic chemical system in which some chemicals are present in high concentrations and some in low concentrations.

Let us consider a chemical systems in which we have three types of chemicals, i.e.  a chemical of type $A $, a chemical of type $B$ and a chemical $K$ that can be at $N$ different states.
We assume that the following chemical reactions take place
\begin{equation} \label{chem system}
\begin{tikzcd}
	{(A)} & {(B)} && {(A)+(k)} & {(B)+(k+1)}
	\arrow["{\rho_0}", shift left=2, from=1-1, to=1-2]
	\arrow["{-\rho_0}", shift left=2, from=1-2, to=1-1]
	\arrow["{\rho_k}", shift left=2, from=1-4, to=1-5]
	\arrow["{-\rho_k}", shift left=2, from=1-5, to=1-4]
\end{tikzcd}
\end{equation} 
Moreover we identify $1$ with $N +1 $.
The state of the chemical system is characterised  by the vector $ (n_A, n_B, k )=(n, k ) \in \mathbb N^2 \times \{ 1, \dots, N \}$.

We assume that the chemical rates $ \{ R_{ \rho_k }(n, k ) , R_{- \rho_k } (n, k )   \}_{ k =0}^N$ at which the reactions occur satisfy the mass action assumption, i.e. we assume that 
\[ 
R_{ \rho_0 }(n, k ) := K_0 n_A, \quad  R_{ - \rho_0 }(n, k ) := K_{-0}  n_B 
\] 
where $K_0, K_{-0} $ are the chemical rates associated with the reactions $ \rho_0, - \rho_0 $ respectively, and 
\[ 
R_{ \rho_k }(n, k ) := K_k n_A, \quad  R_{ - \rho_0 }(n, k ) := K_{-k}  n_B , \quad k \in \{ 1, \dots, N \}.  
\] 
 where $K_k, K_{-k} $ are the chemical rates associtated with the reactions $ \rho_k, - \rho_k $ respectively. 
 
 We associate the stoichiometric vectors $ S_{ \rho_j } \in \mathbb R^{ 2 + N }$ for $ j \in \{ 0, 1, \dots, N \}$ to the chemical reactions $\rho_k $.
We have that 
\[ 
 S_{\rho_0} =- e_1 + e_2 
, \quad S_k = S_{\rho_k} = - e_1 +  e_2 +e_{3+  k}  - e_{ 2+ k}  , \quad k \in \{ 1, \dots, N  \}
\]
and $S_{-\rho_0} =- S_{\rho_0} $ as well as $ S_{ -\rho_k} =  -S_{k}$. 
Here we are denoting with $\mathbf{e}_k $ the $k$-th vector of the canonical basis of $\mathbb R^{2+N}$. 
We want to check if the chemical system \eqref{chem system} satisfies the detailed balance property or not. 
It is easy to notice that the system of reactions \eqref{chem system} contains a cycle $ c =( - N, \textbf{ 1 } )^\dagger   $ where $\textbf{1} \in \mathbb R^N $ is the vector that contains only ones, i.e. $\textbf{1}=(1,1,\dots,1)^\dagger$. In other words the vector  $ c =( - N, \textbf{ 1 } )^\dagger   $ identifies a sequence of reactions that have a net-null effect on the state of the system, indeed
\[ 
\sum_{k=0}^N  c_k S^\dagger _{\rho_k } = \textbf{0}. 
\] 
Therefore the Wegscheider criterion then implies that detailed balance holds only for specific choices of chemical rates (see for instance \cite{feinberg_foundations_2019,wegscheider1902simultane}).
Namely, detailed balance holds if and only if 
\[ 
 \sum_{ j =1 }^N \log\left(\dfrac{K_j }{K_{ - j } }  \right) = N  \log\left(\dfrac{K_0 }{K_{ - 0 } }  \right).  
\] 
Here we assume that the chemical rates $K_{ k }, K_{ - k } $ are such that detailed balance fails. 
 
 The master equation, whose solution $F(t, n , k ) $ is the probability that the system is at state $(n, k )$ at time $ t >0 $, is given by 
 \begin{align*}
     \partial_t F(t, n, k ) & = K_0 ( n_A +1 )  F(t, n + S_0 , k ) - K_{  0 }  n_A F(t, n , k ) +  K_{ -0 }  ( n_B +1 )  F(t, n -  S_0 , k ) \\
     & - K_{ -  0 }  n_B  F(t, n , k )  + K_{k-1} ( n_A +1 )  F(t, n + S_0 , k -1  )  - K_{   k  }  n_A F(t, n , k )\\
     & +  K_{ -(k+1)   }  ( n_B +1 )  F(t, n -  S_0 , k +1 ) -   K_{ -k  }  n_B   F(t, n  , k ) 
 \end{align*}
 where $S_0 =(1, -1 )^\dagger   $. Since $F$ is a probability we also have that 
 \[ 
 \sum_{ n \in \mathbb N^2 } \sum_{k = 1 }^N  F (t , n , k ) = 1 \quad  \forall t \geq 0. 
 \] 
 
 We now assume that $n_A  $ and $n_B $ are of order $ \varepsilon^{-1 } $. More precisely we assume that $ n_A = \dfrac{ c_A }{ \varepsilon } $ and that $ n_B = \dfrac{ c_B }{ \varepsilon } $. 
 In order to obtain non-trivial limits as $ \varepsilon \to 0$ we also assume that $ K_k = \varepsilon a_k $ and that  $ K_{ - k} = \varepsilon b _k $. 
 Finally let us define the function $h$ as 
 \[ 
 h (t, c, k ) := \varepsilon^{-2} F(t, n , k ) 
 \]
 where $c:=(c_A, c_B ) \in \mathbb R_+^2  $.  
 We deduce that the function $h$ satisfies the following equation 
  \begin{align*}
     \partial_t h(t, c, k ) & = \dfrac{ K_0 }{\varepsilon} \left[  ( c_A +\varepsilon  )  h(t, c + \varepsilon S_0 , k ) -   c_A h (t, c , k ) \right]  \\
     & +  \dfrac{ K_{ - 0} }{\varepsilon} \left[ ( c_B +\varepsilon  )  h(t, c -\varepsilon   S_0 , k )   - c_B  h(t, c , k )  \right] \\
     & +   a_{ k-1}  ( c_A +\varepsilon  )  h(t, c + \varepsilon S_0 , k-1  )  -  a_{ k}   c_A h(t, c , k )\\
     & +  b_{k + 1} ( c_B +\varepsilon )  h(t, c - \varepsilon S_0 , k+1 ) -  b_k   c_B   h(t, c  , k ). 
 \end{align*}
 We now use the following approximations valid when $ \varepsilon  $ is small. On one hand we have that 
   \begin{align*}
& \dfrac{ K_0 }{\varepsilon} \left[  ( c_A +\varepsilon  )  h(t, c + \varepsilon S_0 , k ) -   c_A h (t, c , k ) \right]  +  \dfrac{ K_{ - 0} }{\varepsilon} \left[ ( c_B +\varepsilon  )  h(t, c -\varepsilon   S_0 , k )   - c_B  h(t, c , k )  \right] \\
& \approx  S_0 \cdot \nabla_c  [ ( K_0 c_A - K_{-0}  c_B ) h (t, c, k )  ] + \dfrac{\varepsilon }{ 2 }
S_0 \cdot \nabla^2_c  [ ( K_0 c_A +  K_{-0}  c_B ) h (t, c, k )  ] \cdot S_0 \text{ as } \varepsilon \to 0 . \end{align*}
On the other hand 
\begin{align*}
     &  a_{ k-1}  ( c_A +\varepsilon  )  h(t, c + \varepsilon S_0 , k-1  )  -  a_{ k}   c_A h(t, c , k )\\
     &+  b_{k + 1} ( c_B +\varepsilon )  h(t, c - \varepsilon S_0 , k+1 ) -  b_k   c_B   h(t, c  , k ) \\
      & \approx   a_{ k-1}  c_A  h(t, c  , k-1  )  -  a_{ k}   c_A h(t, c , k ) +  b_{k + 1}  c_B h(t, c  , k+1 ) -  b_k   c_B   h(t, c  , k ) \text{ as } \varepsilon \to 0. 
\end{align*}
As a consequence we can write the following approximation for $\partial_t h $ valid as $\varepsilon \to 0 $
  \begin{align*}
     \partial_t h(t, c, k ) & \approx  S_0 \cdot \nabla_c  [ ( K_0 c_A - K_{-0}  c_B ) h (t, c, k )  ] + \dfrac{\varepsilon }{ 2 }
S_0 \cdot \nabla^2_c  [ ( K_0 c_A +  K_{-0}  c_B ) h (t, c, k )  ] \cdot S_0  \\
& +  a_{ k-1}  c_A  h(t, c  , k-1  )  -  a_{ k}   c_A h(t, c , k ) +  b_{k + 1}  c_B h(t, c  , k+1 ) -  b_k   c_B   h(t, c  , k ). 
 \end{align*}
 Since the chemical reactions $ \rho_0 , - \rho_0 $ are much faster compared to the other reactions we can assume that in times of order one the concentrations of $A$ and $B$ approach the stationary concentration $C=(C_A, C_B ) $ satisfying 
 \[
 K_0 C_A - K_{ - 0} C_B =0 \  \text{ and } \ C_A + C_B = C_A(0) + C_B (0)  
 \] 
 where $C_0 =(C_A(0) , C_B (0) ) $ is the initial concentration of $A$ and $B$. 
 Notice that if we isolate the system of macroscopic chemicals $A$ and $B$ then detailed balance always holds.
 Moreover, in times of order one the concentration of $A $ and $B $ will be mainly affected by the faster chemical reaction $S_0$, therefore we can assume that the concentration $c$ of $A$ and $B$ in times of order one is given by 
\[ 
c=c(\xi ) =  C + \xi S_0, 
\] 
where for consistency we must have  
\begin{equation} \label{xi small} 
\vert \xi S_0\vert \ll \vert C\vert.
\end{equation}
We can then change variables and define $f_k (t, \xi ) : = h (t, c (\xi), k ) $. Since 
\[ 
S_0 \cdot \nabla_c = \dfrac{\partial }{ \partial_\xi } \text{ and } S_0 \cdot \nabla_c^2 \cdot S_=  = \dfrac{\partial^2 }{ \partial_{ \xi} } 
\] 
we can write the following equation for $f_k $
 \begin{align*}
     \partial_t f_k (t, \xi  ) & \approx  \dfrac{\partial }{ \partial_\xi }   [ ( K_0 c_A(\xi)  - K_{-0}  c_B(\xi)  ) f_k (t, \xi)  ] + \dfrac{\varepsilon }{ 2 }
 \dfrac{\partial^2 }{ \partial_{ \xi} } [ ( K_0 c_A(\xi)  +  K_{-0}  c_B (\xi) ) f_k (t, c )  ]  \\
& +  a_{ k-1}  c_A(\xi)   f_{ k-1}(t, \xi )  -  a_{ k}   c_A(\xi)     f_{ k}(t, \xi )  +  b_{k + 1}  c_B (\xi)  f_{ k+1}(t, \xi ) -  b_k   c_B(\xi)    f_{ k}(t, \xi )  
 \end{align*} 
  as $\varepsilon \to 0$.
 We can define $\alpha_k(\xi) \coloneqq a_k c_A(\xi) $ and $\beta_k (\xi) \coloneqq b_k c_B (\xi) $.
 Moreover we can define 
 \[
 H'(\xi) := K_0 (C_A + \xi ) - K_{ - 0}  (C_B - \xi ) = (K_0  + K_{ - 0} ) \xi 
 \] 
 and 
 \[
 D(\xi):=  K_0 (C_A + \xi ) + K_{ - 0}  (C_B - \xi ) 
 \]
 Then as $\varepsilon \to 0 $ 
  \begin{align*}
     \partial_t f_k (t, \xi  ) & \approx  \dfrac{\partial }{ \partial_\xi }   [H'(\xi) f_k (t, \xi)  ] + \dfrac{\varepsilon }{ 2 }
 \dfrac{\partial^2 }{ \partial_{ \xi} } [ D(\xi)  f_k (t, \xi )  ]  \\
& +  \alpha_{ k-1}  (\xi)   f_{ k-1}(t, \xi )  -  \alpha_{ k}(\xi)     f_{ k}(t, \xi )  +  \beta_{k + 1}  (\xi)  f_{ k+1}(t, \xi ) -  \beta_k(\xi)    f_{ k}(t, \xi ) . 
 \end{align*} 
 If we assume that $K_0 = K_{ - 0 }  $ then  $H'(\xi)=2 K_0 \xi$ and  $D(\xi)=  K_0 (C_A  +   C_B)  $. 
 
We now expand  $\partial^2_{\xi}(D f_k)=D\partial_{\xi}^2 f_k+  2\partial_\xi D\, \partial_\xi f_k + f_k \partial_{\xi}^2D$. 
 Lastly, since \eqref{xi small} we can write $\delta x =\xi$, for some appropriately chosen $\delta\ll1$.
 Changing coordinates, we see that the dominant term in the diffusion is $\delta^{-2}D\partial^2_{x}f$. To obtain  our model equation \eqref{disc_cont_eq}, we choose $\delta=\varepsilon^{1/4}$ and write $\epsilon=\sqrt{\varepsilon}$.
 
On the other hand, it seems also natural to take $\delta\approx\sqrt{\varepsilon}$, since we expect $c(\xi)$ to be distributed on $C+\R\,S_0$ with variance $\epsilon$ and expectation $0$. Under this assumption, we rewrite the equation as
 \begin{align*}
    \partial_t f_k(t,x)&= \partial_x( (K_0+K_{-0})x\, f_k(t,x))+D(x)\partial_{xx}f_k(t,x)+\sqrt{\epsilon}(2\partial_x D(x)\partial_x f_k(t,x)+\\
     & +[\alpha_{ k-1}  (x)   f_{ k-1}(t, x )  -  \alpha_{ k}(x)     f_{ k}(t, x )  +  \beta_{k + 1}  (x)  f_{ k+1}(t, x ) -  \beta_k(x)  f_{ k}(t,x)].
 \end{align*}
We see that the region $\xi=O(\sqrt{\varepsilon})$, i.e., $x=O(\sqrt{\epsilon})$, forms a boundary layer, in which the diffusion balances the reaction terms.
 
\section{Well-posedness of the time-dependent problem}\label{sec_well_posed}
In this section we prove the existence of solutions to equation \eqref{disc_cont_eq}.
We rewrite the equations for $f$, respectively the relative density $\rho$ with respect to $\gamma$, as follows:
\begin{align}
    \partial_t f &= \partial_x(\epsilon \partial_x f +x\,f)+\partial_x( (H'-x)\,f)+M(x)\,f, \label{full f}\\
    \partial_t \rho &=\epsilon \partial^2_x \rho - x\partial_x \rho +(x-H')\partial_x \rho +M(x)\,\rho .\label{full g}
\end{align}
Let us now define $L_\epsilon \rho=L \rho\coloneqq\epsilon\partial_x^2\rho-x\partial \rho$ and $L^*_\epsilon f=L^* f\coloneqq \partial_x(\epsilon \partial_x+ x f)$.
Let $S^\epsilon_t=S_t$ be the solution semi-group associated to the PDE $(\partial_t-L)u=0$.  This semi-group is known as the Ornstein-Uhlenbeck semi-group, see e.g. \parencite[Section 2]{gaussian}.

It will be convenient to use the notation $a(x)=(x-H'(x))$. 
Then, a formal computation, allows us to write the solution $\rho$ to the full equation \eqref{full g}  by Duhamel's formula as 
\begin{equation} \label{integral rep of P_epsilon}
   \rho(t,x)= P^\epsilon_t ( \rho(t,x)) \coloneqq S_t\,\rho_0(x)+\int_0^t S_{t-s}( a (x) \,\partial_x \rho (t, x)  +M(x)\,\rho (t, x) ) \,ds. 
\end{equation}
The above equality suggests to use a fixed point argument in order to prove the existence of a solution to equation \eqref{full g}. Before proceeding with the fixed point argument, we collect some well-known results for $S_t$, see e.g, {\parencite[Theorem 2.5]{gaussian}}
\begin{theorem}
Let $u(x)\in L^p(d\nu)$, where 
\begin{equation} \label{def: nu} 
\nu_\epsilon=\nu=(1/c)\exp(-x^2/(2\epsilon)).
\end{equation} 
Then we have that 
\begin{equation}
    \vert\vert S^\epsilon_t u\vert\vert_{L^p(d\nu)}\le \vert\vert u \vert\vert_{L^p(d\nu)}
\end{equation}
The operators $S^\epsilon_t \bullet$ are self-adjoint in $L^2(d\nu)$, i.e., for $f,g\in L^2(d\nu)$
\begin{equation}
    \int_\R f\,S^\epsilon_t g \,d\nu=\int_\R g\, S^\epsilon_t f\,d\nu.
\end{equation}
The semi-group has $L=\epsilon\partial_x^2-x\partial_x$ as its infinitesimal generator and has the following integral representation.
\begin{equation}
    S^\epsilon_t u= \int_\R u(y) K^\epsilon_t(x,y)\,dy
\end{equation}
where the integration kernel is given by
\begin{equation*}
    K^\epsilon_t(x,y)=\frac{1}{\sqrt{\epsilon\pi (1-e^{-2t})}}\exp\left(-\frac{\vert y-e^{-t}x\vert^2}{2\epsilon(1-e^{-2t})}\right).
    \end{equation*}
\end{theorem}
\begin{remark}
Since $1\in L^2(d\nu)$, the above in particular shows directly the conservative property of the equation $(\partial_t -L)u=0$.
Since $S^\epsilon_t$ is an operator on $L^p(d\nu)$-spaces it is sometimes convenient to write the kernel with respect to the Gaussian measure $\nu$, that is $K^\epsilon_t (x,y)= N^\epsilon_t(x,y) \sqrt{\pi}^{-1}\exp(-y^2/(2\epsilon))$.
\end{remark}
Now we set up the fixed point argument.
\begin{definition}
Assume that $ u \in L^p(d\nu) $ where $\nu $ is defined as in \eqref{def: nu}. 
Define
\begin{equation*}
    \vert\vert u \vert \vert_{T,\nu, p} \coloneqq \sup_{0\le t\le T} \vert\vert u(t)\vert\vert_{L^p(d\nu)},
\end{equation*}
 The supremum in time should be taken as the essential supremum with respect to the Lebesgue measure. Lastly, we define the Banach space $\C(T,\nu,p)=\{u: ||u||_{T;\nu,p} <\infty\}$. 
\end{definition}
As examples of functions that belong to $\C(T,\nu,p) $ we can consider the functions of the from $u(t,x) = \phi(t)p(x) $, where $\phi(t)$ is smooth and $p(x)$ is a polynomial.
\begin{proposition} \label{prop:existence time dep} 
For every fixed $\epsilon>0$ and $u_0\in L^2(d\nu)$ there exists a (smooth) solution $u \in L^\infty_{loc}(L^2(d\nu))$ of equation \eqref{disc_cont_eq}.
\end{proposition}
\begin{proof}
First, we remove the derivative of $u$ from the above integral representation of $P^\epsilon_t$, i.e. \eqref{integral rep of P_epsilon}. To this end, we integrate by parts, and obtain
\begin{align*}
    S^\epsilon_t (a \partial_x u) &=\int_\R   \partial_y u(y) a(y) K^\epsilon_t (x,y)\,dy \\ 
    &= -\int_\R \left(a(y)\partial_y \log K^\epsilon_t(x,y)+ \partial_y a(y)\right) K^\epsilon_t(x,y) u(y)\,dy. 
\end{align*}
We remark that the term in the second line above is a mapping $L^p(d\nu)\to L^p(d\nu)$. Indeed $\partial_y \log K(x,y) $ is a polynomial and hence belongs to every $L^p(d\nu)$. Thus, for fixed $T\ge 0$, $P^\epsilon_t$ is indeed well-defined from $C(T,\eta,p)$ into itself. Next we estimate $    \vert\vert P^\epsilon_t u_1 - P^\epsilon_t u_2 \vert\vert_{L^p(d\nu)} $ as follows 
\begin{align*}
    \vert\vert P^\epsilon_t u_1 - P^\epsilon_t u_2 \vert\vert_{L^p(d\nu)} &=\vert\vert \int_0^T S^\epsilon_{t-s}( a\,\partial_x u +M\,u)\,ds \vert\vert_{L^p(d\nu)} \\
    &\le \int_0^t \vert\vert S^\epsilon_{t-s} ( a \partial_x u) +S^\epsilon_{t-s}(Mu)\,\vert\vert_{L^p(d\nu)}\,ds\\
    &\le  t\left(\vert\vert M\vert\vert_\infty \sup_{s\le t}\vert\vert u\vert\vert_{L^p(d\nu)}+ C(a) \vert\vert K(.) \vert\vert_{L^q(d\nu)} \sup_{s\le t}\vert\vert u\vert\vert_{L^p(d\nu)}\right)\\
    &\le T\,\max(\vert\vert M\vert\vert_\infty,C(a))\,\vert\vert u \vert\vert_{T,\nu,p},
\end{align*}
where  $1/p+1/q=1$. Therefore, we have a contraction in $\C(T, \nu, p)$ for a sufficiently small $T$, which is allowed to depend on $\epsilon$. Since the equation is autonomous, we can glue local-in-time solutions, and get existence and uniqueness of solutions on $[0,\infty)$, which then are only locally in $L^\infty$. 
\end{proof}

\section{Stationary problem and long-time behaviour}
In this section we prove the existence of stationary solutions to equations \eqref{disc_cont_eq}.  We also prove that the stationary solution in unique. 
Uniqueness of the stationary solutions will be established using the long-time behaviour of solutions of the dual equation, which will be shown to converge exponentially fast to the steady state.
\subsection{Existence}
We start by proving the existence of a stationary solution. 
\begin{lemma}[Existence of invariant measures] \label{lem:existence stationary}
For every initial mass $m>0$, there exists an invariant measure, with density $\Bar{f}=(\Bar{f}_k)_k\geq 0$, $\Bar{f} \in L^1 ( dx )$,
to equation \eqref{full f}, i.e. a solution to 
\begin{equation*}
    \begin{cases}
    0 &=\epsilon\frac{d^2}{dx^2}\Bar{f}+\frac{d}{dx}(H'(x)\bar{f})+M(x)\Bar{f}\\
    m &=\int_\R f\,dx.
    \end{cases}
\end{equation*}
    Moreover, the components of $\bar{f}$ satisfy $\Bar{f}_i \le C_\epsilon \exp(-H(x)/\epsilon)$, where $C_\epsilon=\vert\vert \exp(-H/\epsilon)\vert\vert_{L^1}$.
\end{lemma}
\begin{proof}
For some suitably chosen non-negative initial datum $f_0=g_0\,\nu$, with mass $m$ and $g_0\in L^2(d\nu)$, let $f=(f_k)_k$ be time-dependent a solution to \eqref{disc_cont_eq}. It is easy to see that $f$ preserves the mass and keeps its non-negativity. Thus, we have that $0\le f_k\le \sum_k f_k\eqqcolon F$ and $F$ solves the (scalar) Fokker-Planck equation
\begin{equation} \label{eq:integrated} 
    \partial_t F=\epsilon F_{xx}+(H' F)_x,
\end{equation}
with the constraint that $\int_\R F=m$. This equation has an invariant measure $\gamma\coloneqq Z_\epsilon^{-1}\exp(-H(x)/\epsilon)$ with $Z_\epsilon :=   \int_\R \exp(-H/\epsilon)\,dx $. 
Now, we choose our initial datum $f_0=u\,\gamma$ such that $F_0=U_0\,\gamma$ and that $0 \le U_0 \le C$ is bounded. { Due to Assumption \ref{assumption_H}, we indeed have $g=u (\gamma/\nu) \in L^2(d\nu)$.} By the maximum principle, the solution $U(t)$ to equation \eqref{eq:integrated} with initial datum $U(0)=U_0$ is such that $0\le U(t) \le C$ for all $t$. As a consequence, the family of functions $F(t,\bullet)_{t\ge 0}$ is tight, as $F(t,\bullet)  \le C\gamma(x)$ for every $t \geq 0 $. The same is true for $(f(t,\bullet))_t$. Define now, for $T\ge 0$,
\begin{equation*}
    f_T(x)=\frac{1}{T}\int_0^T f(t,x)\,dt,
\end{equation*}
which is also tight and satisfies $\int_\R f_T \,dx=m$ for all $T$.
Therefore, there exists, up to subsequence, a weak$*$-limit $f_\infty$. To see that $f_\infty$ is a stationary solution we notice that for any test function $\phi$ we have that 
\begin{align*}
   \langle \phi, \LL f_T\rangle&:=\frac{1}{T}\int_0^T \phi\,\left(\epsilon\partial_x^2 f+\partial_x(xf)+M(x)f\right)\,dx\,dt\\
   &=\frac{1}{T}\int_0^T\int_\R \phi(x)\,\partial_t f(t,x)\,dx\,dt =\frac{1}{T}\left(\int_\R \phi\,f(T,x)\,dx-\int_\R\phi\,f(0,x)\,dx \right)\to 0
\end{align*}
as $T \to \infty$ and we take $f_\infty=\bar{f}$. 
The decay estimate follows by construction.
\end{proof}
\begin{corollary}[Existence of $g$] \label{cor:existence g} 
For every $m>0$ and every $\epsilon>0$, there exist non-negative solutions $g=(g_k)_k$ of
\begin{equation}\label{g-eq}
    0 =\epsilon\frac{d^2}{dx} g_k -H'(x)\,\frac{d}{dx} g_k + e_k^\dagger M(x)g,
\end{equation}
that are uniformly bounded in $\epsilon$ and such that $\sum_{k=1}^N g_k(x) \equiv m$.
\end{corollary}
\begin{proof}
Let $f=f_M$ be the stationary solution to \eqref{disc_cont_eq}, provided by the Lemma directly above. Write $f=\gamma\,g$, then $g$ is a stationary solution to \eqref{g-eq}. Since $f\le C_\epsilon \gamma$, the solution $g=g_\epsilon$ is bounded by $C_\epsilon$. Also, $G(x):=\sum_{k=1}^N g_k(x)$ is readily seen to be constant, by looking for bounded solutions of
\begin{equation*}
 \epsilon G_{xx}- H'(x) G_x=0.
\end{equation*}
But if $G(x)\equiv G(0)$ is just a constant function, and $d\gamma$ is a probability measure, then we have $\int_\R F\,dx=G(0)=m$ and $0\le g_k\le G\equiv m$ provides the uniform bounds.
\end{proof}
\subsection{Uniqueness  and spectral gap for \texorpdfstring{$\LL$}{L}}
In order to prove the uniqueness of the invariant measure we consider the \emph{dual equation} to equation \eqref{disc_cont_eq}. Namely we consider the following equation
 \begin{equation}\label{dual_eq}
 \begin{cases}
    -\partial_\tau \phi &=\epsilon \partial_x^2\phi -H'\partial_x\phi +M^\dagger \, \phi=\LL^*\phi,\\
    \phi(T)  &=\phi_0
 \end{cases}
 \end{equation}
where we write $\tau=T-t$.
As it will turn out, the long time behaviour of the solution to the dual equation is easier to understand than the behaviour of the solution to \eqref{disc_cont_eq}, because $\ker M^\dagger (x)=(1,\dots,1)^\dagger $ does not depend on $x$. Along the way, we also establish a spectral gap for the operator $\LL^*$, and exponential convergence to the stationary solutions. We will use this long-time behaviour via the following duality between $f$ and $\phi$.
\begin{align*}
    \int_\R\partial_t(\phi\,f)\,dx=\int_\R \left( \phi \partial_t f(t,x) -f\partial_\tau(T-t,x) \phi \right) \,dx=\int_\R\left(  \phi\LL f- f\LL^*\phi \right) dx=0. 
\end{align*}
If we now integrate from $0$ to $T$ in time, we get 
\begin{equation*}
    \int_\R \phi(0)f(T)\,dx=\int_\R f(0)\phi(T)\,dx.
\end{equation*}
In the end we will send $T\to\infty$ on the right-hand side using the long time behaviour of $\phi$. This allows to obtain information on the long time behaviour of $f$.
\par
We expect that the behaviour of $\phi=\phi(t,x,k)$ in its continuous $x$-variable should be determined by the operator $\epsilon \partial_x^2-H'(x)\partial_x$  while the behaviour in the (discrete) $k$-variables is determined by the matrix $M=M(x)$. This expectation is certainly justified if $M$ does not depend on $x$. Thus, our approach is inspired by the construction of a Lyapunov function for a linear system of the form $\dot{Y}=\tilde{M}\,Y$.
Since in general $M(x)$ is not negative definite (and not symmetric), we would like to have a decomposition of $M(x)$ (resp. $M^\dagger  (x)$).
First, we consider $\tilde{M}$ real and invertible and decompose $\tilde{M}$ as
\begin{equation} \label{eq:splitting} 
    \tilde{M}=S+B, \quad S^\dagger =S, \quad  (SB)^\dagger =-SB,
\end{equation}
with the additional condition $S\in GL(n)$, to exclude e.g. $S=0$. Notice, that $- Y (S/2) Y$ is by construction a Lyapunov function for the system  $\dot{Y}=\tilde{M}\,Y$.
Notice that we can rewrite \eqref{eq:splitting} as the algebraic Riccati equation 
\begin{equation}\label{eq:ricatti}
    0=S^2-(S\tilde{M}+\tilde{M}^\dagger S)/2.
\end{equation}
\begin{lemma} \label{lem: S symmetric and invertible} 
If $\sigma(\tilde{M})\cap-\sigma(\tilde{M})=\emptyset$ there exists a unique matrix $S$ solving the Riccati equation \eqref{eq:ricatti}, and $S$ is symmetric and invertible.
\end{lemma}
\begin{proof}
Consider the related Sylvester equation
\begin{equation}\label{sylvester}
    I= X \tilde{M}+ \tilde{M}^\dagger  X
\end{equation}
which is equivalent to the Riccati equation iff $X=S^{-1}$ exists. If a solution $X$ to \eqref{sylvester} exists, then the symmetric  matrix  $X'=(X+X^\dagger )/2$ is also a solution. Indeed,
\begin{equation*}
    (X+X^\dagger )\tilde{M}+\tilde{M}^\dagger (X+X^\dagger )=(X\tilde{M}+\tilde{M}^\dagger X)+(X^\dagger  \tilde{M}+\tilde{M}^\dagger X^\dagger )=I+I^\dagger =2I.
\end{equation*}
Secondly, if a symmetric solution $X$ exists, then it must be invertible. To see this, choose an orthonormal basis of eigenvectors of $X$. Then, for some orthogonal $Y$, one has
\begin{equation*}
    D N+N^\dagger D\coloneqq D (Y^\dagger  \tilde{M}Y)+ (Y^\dagger  \tilde{M}^\dagger  Y)D=I,
\end{equation*}
where $D=Y^\dagger  \,X\,Y$ is diagonal. If $X$ is not invertible, we may assume that $D\,e_i=0$ for some $i$. Thus,
\begin{equation*}
    e_i=N^\dagger  (D\,e_i)+ D(N\,e_i)=D\,(N\,e_i),
\end{equation*}
which is a contradiction as $e_i$ clearly cannot lie in the range of $D$. Lastly, the spectral condition on $\tilde{M}$, is a well-known existence and uniqueness criterion for solutions of the Sylvester equation.
\end{proof}
Notice that since $ \sigma({M}^\dagger )\cap(-\sigma({M}^\dagger ) )\neq \emptyset $ we cannot apply Lemma \ref{lem: S symmetric and invertible} directly to $\tilde{M} =M(x)$ or to $ \tilde{M} = M^\dagger (x) $. We explain here how we construct a matrix $\tilde{M} $ that satisfies the assumptions of Lemma \ref{lem: S symmetric and invertible}. 

If $S$ is invertible, we can write
\begin{equation*}
    \tilde{M}=S+B=S+ (B\,S^{-1})S,\quad (B\,S^{-1})^\dagger =-(B\,S^{-1}),
\end{equation*}
since 
\begin{equation*}
    B\,S^{-1}= S^{-1}(SB)S^{-1},
\end{equation*}
where $S B$ is skew-symmetric by construction and $S^{-1}$ is symmetric.

\par
 We define 
\begin{align*}
    \ell(x)&\in \ker M(x),\quad s.t.\,\langle \ell(x) , \mathbf{1} \rangle =1
\end{align*}
where we use the notation $(1,\dots,1)^\dagger =\sum_k \mathbf{e}_k \coloneqq \mathbf{1}$ for the generator of $\ker M^\dagger (x)$. We now write the projections onto the kernel and image of $M^\dagger $.  Recall that there is no non-trivial Jordan block for the eigenvalue $0$. Therefore we can define $P $ and $Q $ as follows
\begin{align}
    P(x) &\coloneqq \mathbf{1}\,\ell(x)^\dagger  \label{def P}  \\
    Q(x) &\coloneqq I-P(x) = I-\mathbf{1}\,\ell(x)^\dagger  \label{def Q} 
\end{align}
Notice that both $\ima Q'(x),\,\ima P'(x) \subset \ima P =\ker Q=\ker M^\dagger $. Thus, if we write \eqref{disc_cont_eq} as $\partial_t f= \LL f$ and the adjoint equation as $\partial_t\phi=\LL^*\phi$, the above implies
\begin{equation}\label{Q-LL-relation}
    \partial_t Q\phi=Q\LL^*\phi=Q\LL^* (Q\phi)
\end{equation}
Notice, that by construction $M^\dagger : \ima Q\to \ima Q$ is an isomorphism. Also recall, that all eigenvalues of $M^\dagger \vert_{\ima Q}$ have negative real part and hence the spectral condition stated in Lemma \ref{lem: S symmetric and invertible} to solve the Riccati equation on $\ima Q$ is satisfied. Thus, applying Lemma \ref{lem: S symmetric and invertible}, we deduce that there are projections $\pi(x):\R^n\to\R^{n-1}$ with $\pi\pi^\dagger =I_{n-1}$, $\pi Q=\pi$, such that there exists a matrix $S$ and a matrix $B$ that satisfy 
\begin{equation*}
    \pi M^\dagger  \pi^\dagger  =S+ B,
\end{equation*}
the matrix $S$ is symmetric, invertible  and $L\coloneq BS^{-1}$ is skew-symmetric. Notice that $\pi(x)$ can be expressed in terms of $\ell(x)$, which depends smoothly on $x$ and no component of $\ell(x)$ vanishes anywhere. Thus, we find $\pi(x)$ to be smooth in $x$. 
Essentially, $\pi(x)$ are found by choosing an orthonormal frame for the vector bundle $\ima Q(x)$. Since $\R$ is contractible, any vector bundle is isomorphic to the trivial one, cf. \parencite[Proposition 1.2]{moore2001}.
We can readily see that $S$ is negative definite. 
\begin{lemma}
If $S$ has exactly $n-k$ positive eigenvalues, then $\tilde{M}=S+L S$ has exactly $n-k$ eigenvalues with positive real part.
\end{lemma}
\begin{proof}
For $\delta\in [0,1]$ consider $\tilde{M}_\delta\coloneqq S+ \delta L S$. The matrix $\tilde{M}_0=S$ has inertia $(n-k, k)$ if the number of roots of the characteristic polynomial $p_\delta(x)$ with positive real part of $\tilde{M}_1=\tilde{M}$ is not $n-k$, when counted with multiplicity, then by continuity of the roots there must be a value $\delta$ where there is an eigenvalue with real part $0$. This leads to a contradiction:
Indeed, we find for $\lambda\in\mathbb{C}$ eigenvalue of some $\tilde{M}_\delta$, with eigenvector $0\not=x\in\mathbb{C}^n$ that
\begin{align*}
  \lambda\,x^*S x=(Sx)^*\tilde{M}_\delta x&=x^*(S^2+\delta\, SLS)\,x\\
  &=x^*S^2x+\delta\, x^*SLS x
\end{align*}
Here $x^*(SLS)x\in i\R$, $x^*Sx\in \R$ and $\vert Sx\vert^2>0$, imply that $\lambda$ cannot be purely imaginary, since the right-hand-side has non-zero real part.
\end{proof}
Conjugating $S$ and $L$ by $\pi^\dagger $ and $\pi$, we obtain that 
\begin{equation} \label{A} 
M^\dagger (x) =- \pi^\dagger (x)  A(x) \pi(x)  - \pi^\dagger  (x) L(x)  A(x)  \pi(x)
\end{equation} 
where we set $A (x)\coloneqq -S(x)  \ge 0$. This is the matrix decomposition we will use in the following.
\begin{proposition}[Long-time limit of the dual solution] \label{long_time_dual}
Denote $\gamma(x)= c_\epsilon \exp(-H/\epsilon)$ for every $ x \in \mathbb R$.
Let $\phi(t,\bullet)=\phi(T-t,\bullet )\in L^2(\gamma\,dx)=:L^2(d\mu)$ be a solution to the dual equation. Then, there exist $\epsilon_0>0$, $a_\infty >0$ and $\lambda>0$, such that for each $\epsilon\le \epsilon_0$,
\begin{equation*}
    \vert\vert \phi(t)- a_\infty \mathbf{1}\vert\vert_{L^2(d\mu)}\le C \exp(-\lambda t).
\end{equation*}
If $H$ is uniformly convex, $\lambda$ can be taken independent of $\epsilon$.
\end{proposition}
\begin{proof}
 We write $\phi_Q(t, x) =Q(x)\phi(t,x)$ and $a(t,x)= \phi(t,x)\cdot \ell(x)$. By the definition of $P(x)$ and $Q(x)$ we have that $P(x)\phi(t, x ) =a(t,x)\mathbf{1}$ as well as $\phi=\phi_Q+a \mathbf{1}$. As 
observed in \eqref{Q-LL-relation}, we have $Q\LL^*=Q\LL^*Q$. This yields
\begin{align*}
\int_\R \langle A\phi_Q, Q\LL^*\,Q\phi_Q\rangle\,d\mu &=\int_\R \langle A\phi_Q, \partial_x(\epsilon\partial_x\phi_Q\,\gamma)+M^\dagger \phi_Q\rangle\,dx\\
    &=-\int_\R \epsilon\left(\langle A\partial_x\phi_Q,\partial_x\phi_Q\rangle + \langle \phi_Q \partial_x A ,\partial_x\phi_Q\rangle\right)\,d\mu-\int_\R \langle A^2 \phi_Q,\phi_Q\rangle\,d\mu 
\end{align*}
where $A(x)$ is given by \eqref{A}. 
The only summand which is not already non-positive is the one involving $A'(x)$.
Using Cauchy-Schwarz and Young's inequality, we get
\begin{equation*}
    \vert \langle A'(x) \phi_Q(t, x), \partial_x\phi_Q(t, x)\rangle\vert\le \dfrac{1}{\delta} \vert\vert A'\vert\vert_\infty^2\vert \phi_Q(t, x)\vert^2 + \delta\vert \partial_x \phi_Q(t, x)\vert^2, \quad \delta >0 .
\end{equation*}
We choose $\delta>0$ small enough to ensure $v^\dagger  A(x) v - \delta \vert v\vert^2> 0$ for $v\in\ima Q(x)$ (i.e., $\delta$ smaller than the minimal spectral gap of $S(x)$). Lastly, we take $\epsilon$ small enough to ensure that
\begin{equation*}
    \frac{\epsilon}{\delta}\vert\vert S'\vert\vert_\infty I - S^2(x)< 0,
\end{equation*}
as bilinear forms. 
These choices of $\epsilon>0 $ and $\delta>0 $ guarantee that 
\begin{align*}
    \epsilon \vert \langle \phi_Q \partial_x A ,\partial_x\phi_Q\rangle \vert \leq \dfrac{\epsilon }{\delta } \vert\vert A' \vert\vert_\infty   \vert \phi_Q \vert^2+ \epsilon \delta \vert \partial_x \phi_Q \vert^2 <  \langle A^2 \phi_Q,\phi_Q\rangle + \epsilon  \langle A\partial_x\phi_Q,\partial_x\phi_Q\rangle. 
\end{align*}
Recall that $A$ is positive definite on $\ima Q$, and that $\langle v, A(x)v\rangle^{1/2}$ is equivalent to the euclidean norm on $\R^n$. Therefore, there exists $s_Q>0$ such that 
\[
\int_\R \phi_Q\cdot Q\LL^*\phi_Q\,d\mu\le -  s_Q \int_\R \phi_Q^2\,d\mu
\]
and by Gr\"onwall
\[\vert\vert \phi_Q(t)\vert\vert_{L^2(d\mu)}^2\le \vert\vert \phi_Q(0)\vert\vert_{L^2(d\mu)}^2\exp(-s_Q t).\]
Treating $\phi-\phi_Q$ is more involved because $P\LL^* \phi\not = P\LL^*P\phi$, as $\partial_x Q(x)\in \ima P(x)$. Therefore, $a(t,x)= \phi(t,x)\cdot \ell(x)$ does not solve the scalar equation $(\partial_t-\epsilon\partial_{xx}+H'\,\partial_x) a =0$.
Thus, we write \eqref{dual_eq} as an inhomogeneous equation susing $\phi_Q$ and then employ  Duhamel's formula. We write $M^\dagger (x)\phi=M^\dagger  \phi_Q\eqqcolon R(t,x)$, where we just proved that $\vert\vert R(t)\vert\vert_{L^2(d\mu)}^2=O(\exp(-s_Q t))$ as $t \to \infty $. Thus, $\phi$ solves the inhomogeneous equation
\begin{equation*}
    (\partial_t-\epsilon\partial_{xx}+H'\partial_x)\phi(t, x) =R(t,x)
\end{equation*}
and we write for $T_0\ge 0$
\begin{equation}\label{duhamel_phi}
    \begin{aligned}
    \phi(T_0+t,x)&=P^\epsilon_t \phi(T_0, x ) + \int_{T_0}^{t+T_0} P_{t+T_0-s}^\epsilon R(s,x )\,ds\\
    &=P^\epsilon_t a(T_0,x)\mathbf{1} +P^\epsilon_t \phi_Q(T_0, x ) + \int_{T_0}^{t+T_0} P_{t+T_0-s}^\epsilon R(s, x )\,ds,
    \end{aligned}
\end{equation}
where the semigroup $P^\epsilon_t$ is generated by $\epsilon\partial_{xx}- H'\,\partial_x$.
Observe that $\vert\vert P_t^\epsilon g\vert\vert_{L^2(\mu)}\le \vert\vert g\vert\vert_{L^2(\mu)}$ for all $g\in L^2(d\mu)$ and that by the definition of $P_t\epsilon $  we have that for every $g \in L^2 (d \mu ) $
\begin{equation} \label{invariance P_t}
\int_\R P_t g \,d\mu=\int_\R g \,d\mu. 
\end{equation} 
In particular, $\phi(t)=P_t^\epsilon \phi(0) + \int_{0}^t P_{t-s}^\epsilon R(s)$ is bounded in $L^2(d\mu)$ and we find
\begin{equation*}
    \vert\vert \phi(T_0+t)-P_t^\epsilon (a(T_0)\mathbf{1})\vert\vert_{L^2(d\mu)} \le C ||\phi_Q(0)||_{L^2(d\mu)}\exp(-s_Q T_0)
\end{equation*}
for a suitable constant $C>0 $. 
We may rewrite this as 
\begin{equation*}
    \vert\vert a(t+T_0,x) - P_t^\epsilon a(T_0)\vert\vert_{L^2(d\mu)} \le C\exp(-s_Q T_0).
\end{equation*}
Now, we can use the decay of $P_t^\epsilon a(T_0)$ as $t\to\infty$.
If $H(x)$ is (uniformly) convex on $\R$, we can use a well-known Poincaré inequality due to \cite{bakryemry}, \parencite[Theorem 3.1]{schlichting_menz}. There exists $C_P$ such that
\[\Var_\mu(f)\le \epsilon C_P\,\int f'^2\,d\mu,\]
with $\Var_\mu(f)=\int (f-(\int f\,d\mu))^2\,d\mu$, and $C>0$ independent of $\epsilon$. By linearity we can replace $P_t^\epsilon a(T_0)$ by $\bar A(t)=P_t^\epsilon a(T_0)-\int_\R P_t^\epsilon a(T_0)\,d\mu$ and obtain 
\begin{equation*}
       \int_\R \bar A\partial_t  \bar A \,d\mu=  \int_\R \bar A(\epsilon \partial_x^2 \bar A -H'\partial_x \bar  A)\,d\mu= \epsilon \int_\R \bar A \, \partial_x( \gamma\partial_x \bar A )\,dx
       = - \epsilon \int_\R \bar (\partial_x \bar A )^2 \,d\mu \le -C_P \int_\R \bar A^2\,d\mu.
\end{equation*}
We deduce that 
\begin{equation} \label{bound for bar A}
    \| \bar A (t)  \|_{L^2(\mu)} \leq C \exp \left( - C_P t \right) 
\end{equation}
for $C>0 $. 
For any $f$ and $t\ge 0$, we have $\Var_\mu(P_t^\epsilon  f)\equiv \Var_\mu(f)$.
 We finish the proof by Gr\"onwall since we can pass to the limit $\lim_{T_0\to\infty}\int_\R a(T_0)\,d\mu$, by
integrating \eqref{duhamel_phi} for some fixed time $T=T_0+t$.
{
Indeed, we have that
\begin{align*}
  a_\infty &\coloneqq \lim_{T\to\infty} \int_\R \phi(T)\,d\mu=\lim_{T\to\infty} \int_\R P_T \phi(0)\,d\mu +\lim_{T\to\infty} \int_\R \int_0^T P_{T-s} R(s) \,d s \, d\mu\\
  &= \int_\R \phi(0)\,d\mu + \int_0^\infty \int_\R P_{T-s} R(s) \,d\mu \, ds 
\end{align*}
converges, because $\int_\R R\,d\mu$ is integrable in time and $P^\epsilon_\bullet$ preserves the mass. In the second line of \eqref{duhamel_phi}, we have
\begin{equation*} 
  \int_\R \phi(T_0+t)\,d\mu -\int_\R P_t a(T_0)\,d\mu=\int_\R P^\epsilon_t \phi_Q(T_0, x ) + \int_\R \int_{T_0}^{t+T_0} P_{t+T_0-s}^\epsilon R(s, x )\,ds= O(\exp(-s_Q T_0),
\end{equation*}
as $T_0\to \infty$.
 Therefore, for all $t\ge 0$, $a_\infty=\lim_{T_0\to\infty} \int_\R P_t a(T_0)\,d\mu$. Lastly, we recall that $\int_\R P_t a(T_0)\,d\mu=\int_\R  a(T_0)\,d\mu$.
}
Using \eqref{duhamel_phi}, as well as the definition of $a$, the facts that $\vert\vert R(t)\vert\vert_{L^2(d\mu)}^2=O(\exp(-s_Q t))$ and that $P_t^\epsilon$ is a contraction  
\begin{align*} 
    \vert\vert \phi(t)- a_\infty \mathbf{1}\vert\vert_{L^2(\mu)} & \leq  \vert\vert  P_{t/2}^\epsilon \phi(t/2)- a_\infty \mathbf{1}\vert\vert_{L^2(\mu)} +
    \int_{ t /2 }^t  \vert\vert P_{ t -s }^\epsilon R(s) \vert\vert_{L^2(\mu)}  ds \\
    & \leq \vert\vert  ( P_{t/2}^\epsilon a(t/2)- a_\infty ) \mathbf{1}\vert\vert_{L^2(\mu)} +
    C_1 \exp \left(- s_Q t/2  \right).  \\
\end{align*} 
In order to estimate $\vert\vert  ( P_{t/2}^\epsilon a(t/2)- a_\infty ) \mathbf{1}\vert\vert_{L^2(\mu)}$ we use \eqref{bound for bar A} and \eqref{invariance P_t} 
\begin{align*} 
\vert\vert  ( P_{t/2}^\epsilon a(t/2)- a_\infty ) \mathbf{1}\vert\vert_{L^2(\mu)} & \leq \vert\vert  \bar A (t/2)  \mathbf{1}\vert\vert_{L^2(\mu)}  + \vert\vert   \int_\R ( P_{t/2}^\epsilon a (t/2) d\mu - a_\infty ) \mathbf{1}\vert\vert_{L^2(\mu)} \\
& \leq C_2 \exp \left( - C_P t/2 \right) +  \vert\vert  ( \int_\R a (t/2) d\mu - a_\infty ) \mathbf{1}\vert\vert_{L^2(\mu)}\\
&\le  {C_2 \exp \left( - C_P t/2 \right)+C_3 \exp(-s_Q t/2)}
\end{align*}

Combining the two estimates shows the claim
\begin{equation*}
    \vert\vert \phi(t)- a_\infty \mathbf{1}\vert\vert_{L^2(\mu)} \le (C_1+C_3) \exp(-s_Q t/2) + C_2 \exp(-C_P t/2) \le C \exp(-\lambda t),
\end{equation*}
where $\lambda\coloneqq \min(s_Q, C_p)/2$.
\end{proof}
\begin{remark}
If $H$ is allowed to be non-convex, e.g., a multi-well potential, the spectral gap becomes exponentially small in $\epsilon$. If $H$ is not convex, we can employ \parencite[Corollary 2.15]{schlichting_menz}, which provides exponentially large (in $\epsilon$) Poincaré constant $C_P(\epsilon)$. In this case, the conditions outlined in Remark \ref{rem:assumption_H} is sufficient to apply these results.
\end{remark}
As a corollary of the exponential return to the stationary solution, we find a uniform spectral gap for the operator $\LL^*$.
\begin{theorem}[Spectral Gap for $\LL^*$]
We have $\sigma_p(\LL^*)\subset \{0\} \cup[-\lambda,-\infty)$, where $\lambda$ is the same as in  Proposition \ref{long_time_dual}
\end{theorem}
\begin{proof}
Suppose there is an eigenfunction of $\LL^*$ $\psi\in L^2(d\mu)$ for $\eta> -\lambda$, then $\exp(\eta t)\psi(x)$ is a solution to equation \eqref{dual_eq}. However, it is easy to see that Proposition \ref{long_time_dual} is violated for any candidate $a_\infty\in\R$, a contradiction.
\end{proof}
\begin{remark}
It is not the case that we can always find a basis of eigenvectors for $L^2(d\mu)$, even for constant matrix $M$. Recall $M$ is in general neither symmetric nor diagonalisable.
\end{remark}
Lastly, we use the convergence of the dual solution to establish the uniqueness of the stationary solutions of \eqref{disc_cont_eq}.
\begin{lemma}[Uniqueness of the invariant measure]
 For every $m>0$, the stationary solution $\bar{f}$ of \eqref{disc_cont_eq} with mass $m$ and $\bar{g}\in L^2(d\mu)$ is unique.
\end{lemma}
 \begin{proof}
 Fix $T>0$, then we have, by integration by parts,
 \begin{align*}
    \int_\R f(T,x )\phi(0,x )\,dx - \int_\R \phi(T,x )f(0,x )\,dx &=\int_0^T\int_\R \partial_t(f(t,x)\phi(t ,x ))\,dx \,dt\\
    &=\int_0^T\int_\R \left( f(t, x) \LL^*\phi(t, x) - \phi(t, x)\LL  f(t, x)  \right)  \,dx \,dt\\
    & =0.
 \end{align*}
 Assume that $h=f_1-f_2$ where $f_i=\gamma\,g_i$, with $0\le g_i\le C$, are two stationary solutions with the same mass $m>0$.  Let $\phi(0)$ be any (vector-valued) compactly supported test function. Then $\phi(0)\in L^2(d\mu)$ and we get
 \begin{equation*}
     \int_\R h (x)\phi(0,x)\,dx=\lim_{T\to\infty}\int_\R h (x) \phi(T,x)\,dx=a_\infty\int_\R (f_1(x) -f_2(x) )\cdot \mathbf{1}\,dx=0.
 \end{equation*}
 Thus, we have $f_1=f_2$ and the stationary solution $\bar f $ is unique.
 \end{proof}
  \begin{remark}
 Since every bounded solution $g_\epsilon$ defines an invariant measure for \eqref{disc_cont_eq},  writing $f\coloneqq g_\epsilon(x)\,\gamma$ immediately also shows uniqueness of the solution $g$ to equation \eqref{g-eq}.
 \end{remark}
 Essentially the same proof also shows the exponential convergence of the solution $\rho(t,x)$  to equation \eqref{disc_cont_eq} with initial mass $m>0$ to the (unique) stationary solution with the same mass $m>0$.
 \begin{corollary}
 Let $\rho \in L^\infty_{loc}(L^2(d\mu))$ be a solution to \eqref{disc_cont_eq rho}. Let further $\bar{f}=g \gamma$ be the stationary solution with mass $m$. Then
 \begin{equation*}
     \vert\vert \rho(t)-g\vert\vert_{L^2(d\mu)}\le C\exp(-\lambda t).
 \end{equation*}
 Recall that $\partial_t \rho= (\gamma^{-1} \LL \gamma)\rho$. Therefore, in $L^2(d\mu)$, we also have $\sigma_p(\gamma^{-1}\LL \gamma)\subset \{0\} \cup [-\lambda,-\infty)$
 \end{corollary}
 This Corollary is Theorem \ref{thm:main reults spectral gap}. 
 \section{Preliminary results on \texorpdfstring{$M(x)$}{M(x)}} \label{sec_M}
To understand the dynamics of the $f_k$, we need some observations on the structure of the discrete part of the system.
The following are simple consequences of the well-known Perron-Frobenius theorem, cf. \parencite{SpectralProofPerronFrobenius, perron1907}, and some general observations about the form of the matrix $M(x) $. 
\begin{lemma}\label{lemma_M}\hfill

Assume that $M:\R\to \R^{n\times n}$ satisfies Assumption \ref{assumption_M_structure}. Then for every $ x \in \mathbb R$ we have that 
\begin{enumerate}[label=(\roman*)]
    \item\label{simpleEV}  $0$ is a simple eigenvalue of $M(x)$, i.e., there are no non-trivial Jordan blocks for $0$.
    \item $\mathbf{1}\coloneqq \sum_k \mathbf{e}_k\in \ker M(x)^\dagger $, hence no eigenvalue of $M(x)$ except $0$ can have (component-wise) a non-negative eigenvector.
    \item If $\lambda \not =0$ is an eigenvalue of $M(x)$, then $\Re \lambda <0$.
    \item\label{item_maximum_principle} If all $v =(v)_{i=1}^n \in [0, \infty)^n  $ and $j \in \{ 1, \dots, n \} $ is such that $v_j=0$ then $(M(x) v)_j\ge 0$.
    \item \label{no cycles}There is $\ell(x)\in \ker M(x)$ such that $\ell(x)_i > 0$ for every $ i \in  \{ 1, \dots, n \} $.
\end{enumerate}
\end{lemma}
\begin{proof}
Keep $x$ fixed and write $M(x)=M$. Most of the statements follow from Perron-Frobenius and the observation $m_{i,j}\ge 0$ whenever $i\not=j$. 
Point \ref{simpleEV} requires the most work. Let $v\in \R^n$ be a generalised eigenvector such that $M^r v=0$ and $M^{r-1}v\not=0$, then by adding some $\ell\in \ker M$, we may assume that all $v_k\ge 0$. Consider the ODE $\dot{X}=MX$, with initial condition $X(0)=v$. Then, by \ref{item_maximum_principle}, we see $X(t)\ge 0$ component-wise for all $t>0$.
On the other hand, since $(1,\dots,1)^\dagger \in \ker M^\dagger $, we have that $0\le X_k(t)\le \sum_k X_k(t)\equiv \sum X_k(0)$. In particular the solution $X(t)$ is bounded. We have $X(t)=\exp(t M)v$ and this is readily seen to be a polynomial of degree $r-1$ in $t$. However, such a polynomial can only be bounded for $r=1$.
\end{proof}
\begin{remark}
In view, of the above proof, Point \ref{item_maximum_principle} should be understood as a kind of discrete maximum principle (or minimum principle) in the variable $k$. If the non-negative function $v:\{1,\dots, N\}\to \R_+$ assumes a local minimum $v_k=v(k)=0$, then, the condition $(M(x)v)(k)\ge 0$ replaces the non-negativity of the second derivative at local minima of differentiable functions. From this we deduce, e.g, that $f_k(0,x)\ge0$ implies positivity of the solutions to Equation \eqref{disc_cont_eq} for all $t>0$.
\end{remark}

 \section{Asymptotic behaviour of the limit equation}\label{sec_asymptotics}
Away from the points where $H'=0$, at which we will have to introduce a boundary layer, we expect the asymptotics of $g_\epsilon$ as $\epsilon \to 0 $ to be given by $ g $ where $g$ is the solution of the limit equation
\begin{equation}\label{limit_eq}
    H'(x)\,\frac{d}{dx}g(x)=M(x)\,g (x).
\end{equation}
The existence of a (bounded) global solution to the above limit equation is not an immediate consequence of Picard-Lindelöff Theorem, since $ x \mapsto (1/H'(x) )M(x)$ is not Lipschitz at the critical points.
In this section we study the existence of a solution to \eqref{limit_eq} and its asymptotic behaviour as $x \to 0$. 
We consider the intervals $ (0, \infty ) $ and $ ( - \infty , 0 )$ separately.
On the interval $(0,\infty)$ it is convenient to rewrite equation \eqref{limit_eq} using a change of variables. 
In particular,  we consider the change of variables $y(x)$ such that 
\begin{equation*}
    \dfrac{d y(x)}{dx}= \dfrac{1}{H'(x) }. 
\end{equation*}

For example, if $H(x)=\dfrac{x^2}{2}$, then $ y(x)= \log(x)$. For generic $H$ satisfying Assumption \ref{assumption_H}, we have a similar behaviour. 
Since $H''(0)>0$, we have that $x\sim \exp(c y)$ for a constant $c >0$ as $y\to-\infty$ while for $y\to +\infty$ we use our assumption $x-C\le H' \le x+C$ to obtain the same behaviour.
 Therefore, we consider now the following equation
 \begin{equation*}
     \frac{d}{dy} \bar g(y) = \bar{M}(y)\bar{g}(y),\quad \lim_{y\to-\infty} \bar g(y)=g_0.
 \end{equation*}
  where $ \bar M (y):= M (y (x) )$ and $ \bar g (y):= g (y(x)) $. $\bar{P}$ and $\bar{Q}$ are defined analogously. In the following we will omit the bar, if no confusion can arise. We write
  \begin{align*}
      M_-&\coloneqq\lim_{y\to-\infty} \bar M(y)=M(0),\,\text{and}\\ M_+&\coloneqq\lim_{x\to\infty}M(x)=\lim_{y\to \infty} \bar{M}(y);
  \end{align*}
   likewise define $P_\pm$ and $Q_\pm$ for the projections $P(x)$ and $Q(x)$ in  \eqref{def P} and \eqref{def Q}. 
 
$M(y)$ is  bounded and smooth, hence solutions $g(y)$ with any initial datum $g(y=0)$ exist on the whole of $\R$, but may diverge at $-\infty$, because the eigenvalues $M(y)$ only have negative real part. 
Corollary \ref{cor:existence g} also implies that the stationary solutions of \eqref{g-eq} satisfy
\begin{align*}
    \frac{d}{dy}\left(\sum_{k=1}^n g_k (y) \right)=0.
\end{align*}
By comparison and point \ref{item_maximum_principle} of Lemma \ref{lemma_M} we know that if $g(y=y_0)_k\ge0$, then the same holds for all $y>y_0$ and solutions remain bounded.
The result of this section is the following statement.
\begin{proposition} \label{prop: existence and uniqueness solution limiting problem} 
For each $m>0$, there exists a unique global solution $g=(g_k)_k \ge 0$ to
\begin{equation*}
    \frac{d}{dy} g(y) = M(y) g(y),\quad y \in \R 
\end{equation*}
such that for all $k=1,\dots,n$
\begin{align*}
     \sup_{y\in \R} \vert g_k(y) \vert \le m.
\end{align*}
Further, 
\begin{equation*}
    \lim_{y\to-\infty} g(y) \in \ker M_-  
\end{equation*}
where $M_-:= M(x=0)$.
Moreover, for every (potentially unbounded) non-negative solution of this equation the limit  
\begin{equation*}
    \lim_{y\to \infty} g(y) \in \ker M_+, 
\end{equation*}
exists, where $M_+:= \lim_{ y \to \infty } M(y)$. 
Both limits are uniquely determined in $\ker M_\pm$ by the constant $m=\sum_k g_k(x)$.
\end{proposition}
\begin{proof}
We begin by constructing the unique bounded solution.
Since we are interested in studying the behaviour of $ g $ as $ y \to - \infty $, it is convenient to use the change coordinates to $y=-z$. Then we have that $\tilde{g} (z)=g(-z) $ satisfies $ \dfrac{d \tilde{g}(z) }{d z } =- M(-z ) \tilde{ g } (z) $. We define the matrices 
\begin{equation*}
A := -M_- \ \text{ and } \  B(z):= -(M(-z )-M_-). 
\end{equation*}
Then we have that 
\begin{equation} \label{ODE g tilde} 
\dfrac{d \tilde{g}(z) }{d z } = (A + B(z)) \tilde{g} (z) . 
\end{equation} 
Notice that by definition $A$ has now only eigenvalues with positive real part on $\ima Q_-$.
Moreover notice that $B(z)=O(\exp(-z) )$ as $ z\to \infty$.  In order to see this we can first of all observe that if we have that 
\begin{equation*}
\left.\frac{d}{dx}M(x)\right\vert_{x=0}=M'(0)\not=0,
\end{equation*}
then $M(x)\approx M_-+ M'(0) x+O(x^2)$ as $x\to0$. Which in turn implies the desired result.
If $M'(0)=0 $ one can obtain the same result considering higher orders in the Taylor expansion.

Now consider, $a\in C^0_b(\R_+)$, $z_0>0$ and $ g_0\in \ker A=\ima P=\ker Q$. Let us define the functional $I : C^0_b(\R_+) \rightarrow C^0_b(\R_+) $ obtained applying Duhamel's formula to \eqref{ODE g tilde}
\begin{equation*}
    I[a](z)=I[a,g_0,z_0](z)\coloneqq  e^{ (z - z_0) A } g_0  + \int_{z_0}^z e^{A(z-s)} P\,B(s)\,a(s)\,ds - \int_{z}^\infty e^{A(z-s)} Q \,B(s)\,a(s)\,ds.
\end{equation*}
A fixed point $I(g,z_0,g_0)=g$ will be a solution of the equation $d/dz g= M(z) g(z)$, with initial datum $g(z_0=g_0+\int_{z_0}^\infty  e^{A(z-s)} Q \,B(s)\,a(s)\,ds$, i.e., $Pg(z_0)=g_0$.
Notice that $\exp(t A)P=P$ for all $t$, since $\exp(t A)\vert_{\ker A}=I$. Therefore, $e^{ (z - z_0) A } g_0=g_0$ and the first integrand $e^{A(z-s)} P\,B(s)\,a(s)$ is uniformly integrable for bounded $a$ and choosing $z_0$ large enough one can bound the integral by, say $1/2\,\vert\vert a\vert\vert_\infty$, since $B(s)=O(\exp(s))$.\par
Similarly, for all $t\in \R_+$, $\exp(-t A)\,Q$ is contracting at least like some $t^N\exp(-c t)I$ and the latter integral is always finite. Thus, we can make the latter integral smaller that $1/2 \vert\vert a\vert\vert_\infty$ by choosing $z_0\le z$ large enough and using the decay of $B(z)$.\par
Therefore, using the fact that $I[\bullet]$ is a contraction, we conclude that for any $g_0\in \ker M_-$ there exists a unique solution $g'(y)= M(y) g(y)$ with $P g(-z_0)=g_0$, and that this solution remains bounded for $y\to -\infty$. It is straightforward to see, using the fixed point equation $I[g]=g$, that $g_*\coloneqq\lim_{y\to-\infty} g(y) \in \ker M_-$ must exist.  Since  $\ker M_-$ is one-dimensional, $g_*$ is uniquely determined by $\langle g_*,\mathbf{1}\rangle \eqqcolon m>0$. We recall that
\begin{equation*}
    \frac{d}{dz}\sum_k g_k(z) =0,
\end{equation*}
and also that $\langle \mathbf{1}, Q v\rangle =0$. Thus, $m= \langle \mathbf{1}, P g(z)\rangle\equiv \langle \mathbf{1}, g_0\rangle$, and hence $g_*=g_0$.
The non-negativity of the $g_k$ follows from Lemma \ref{lemma_M}, as elements in $\ker M(x)$ are component-wise positive for all $x\in \R$.

Now we consider the limit $y\to\infty$. For this, we use a similar linearisation of $M(x)$ as in the above case of $y\to -\infty$.
Take any $y_0 >0 $. 
We define $E(y)= \sum_{k=1}^n g_k(y)\equiv E(y_0)$. Notice that if $g_k(y_0)\ge 0$, then we have also $0\le g_k\le E(0)$ for all $y\in[y_0,\infty)$.
As before it is convenient to write $M(y)=M_+ +(M(y)-M_+)=M_+ +R(y)$. For large $y>0$, one can write a variation of constants formula for the inhomogeneous ODE $\frac{d}{dy} G = M_+ G+ R(y)g$
\begin{equation*}
    g(y)= \exp(M_+(y/2))g(y/2) + \int_{y/2}^y \exp(M_+\,(y-s) )\, R(s)g(s)\,ds. 
\end{equation*}
Using the boundedness of $R(s)$ and $g(s)$ together with the fact that $M_+$ has negative eigenvalues on $\ima Q_+$ with $Q_+ = \lim_{ x \to \infty } Q(x) $, we see that the integral above is bounded, up to a multiplicative constants, by $\max_{y>\overline{y}} \vert\vert R(y)\vert\vert$, which vanishes as $\overline{y}\to\infty$. In particular, we have
\begin{equation*}
    \lim_{y\to\infty} g(y) =g_+ \in \ker M_+,\quad E(\infty)\coloneqq \sum_k \langle \mathbf{e}_k, g_+\rangle =E(y_0),
\end{equation*} 
since $E(y_0)$ determines a unique element in $\ker M_+$.
\end{proof}
\begin{remark}
Using a completely analogous change of variable, we can repeat this procedure for $x\in (-\infty,0]$ to obtain a solution to \eqref{limit_eq} that is defined for $x$ in the whole of $\R$.
\end{remark}
\begin{remark}
For general, non-convex, Morse function $H$ it is possible to consider intervals of the form $(p_i,q_i)$ and $(q_i, p_{i+1})$, where $p_i$ are the local minima of $H$ and $q_j$ the local maxima in between the minima.
We can then repeat the arguments used in this section on these intervals and deduce analogous results on each of the intervals.
\end{remark}
\section{Convergence to the limit equation}\label{sec_convergence}
This section is concerned with proving the  convergence of the stationary solutions $g_\epsilon$ of \eqref{g-eq} to the unique solution of the limit equation \eqref{limit_eq} constructed in Section \ref{sec_asymptotics}.

For our purposes it is convenient to write equation \eqref{g-eq} in an alternative way. 
We define $h_\epsilon \coloneqq d/dx { g_\epsilon}$. Then \eqref{g-eq} can be written as an inhomogeneous system
\begin{align*}
    \frac{d}{dx}  g_\epsilon  &= h_\epsilon(x) \quad x \in \R  ,\\
    \frac{d}{dx}  h_\epsilon (x) &= \frac{H'(x)}{\epsilon}h_\epsilon (x)  + \frac{1}{\epsilon} F(x) ,\quad F(x)\coloneqq M(x)\,g_\epsilon (x).
\end{align*}
Therefore, we may use  variation of constants formula to obtain that 
\begin{equation*} 
    h_\epsilon(x)=\exp(H(x)/\epsilon)\left(h_0 +\int_0^x \frac{1}{\epsilon}\exp(-H(s)/\epsilon)\,F(s)\,ds \right), \quad  x >0,
\end{equation*}
where we are taking $\epsilon>0$ sufficiently small, but fixed. 
However notice that we need to impose a compatibility condition for $h_0\coloneqq h_\epsilon(0)$. The reason is that we are looking for bounded solutions $g$.
Since $M(s)$ and $g(s)$ are bounded, so is $F(s)$. In particular, 
\begin{equation*}
 h_\infty\coloneqq\int_0^{\infty}\frac{1}{\epsilon} \exp(-H(s)/\epsilon)\,F(s)\,ds\,=C(F,\epsilon)<\infty.
\end{equation*}
This allows to rewrite the variation of constants formula for $h_\epsilon $ as follows
\begin{equation*}
    h_\epsilon(x)=\exp(H(x)/\epsilon)\left( (h_0+h_\infty)-\int_x^\infty \frac{1}{\epsilon} \exp(-H(s)/\epsilon)\,F(s)\,ds\, \right)
\end{equation*}
It is easy to see that the term $e^{H(x)/\epsilon}\int_x^\infty \epsilon^{-1}e^{-H(s)\epsilon}F(s)\,ds$ remains bounded as $x\to\infty$. Thus, in order to have that  $h$ and $g$ are bounded we need to have $h_0=-h_\infty$.
This compatibility condition allows us to reduce the equation for $h_\epsilon$ to 
\begin{equation}\label{eq:int for h_epsilon}
    h_\epsilon(x)= - \exp(H(x)/\epsilon)\int_x^\infty \frac{1}{\epsilon} \exp(-H(s)/\epsilon)\,F(s)\,ds\, 
\end{equation}

The analogous argument applies for the integral from $-\infty$ to $0$, giving $h_0=h_{-\infty}$. Indeed, for $-x$ with $x>0$, we find,
\begin{align*}
\int_0^{-x} \frac{1}{\epsilon} \exp(-H(s)/\epsilon)\,F(s)\,ds &=-\int_{-x}^0 \frac{1}{\epsilon} \exp(-H(s)/\epsilon)\,F(s)\,ds\\
&=-h_{-\infty}+\int_{-\infty}^{-x}\frac{1}{\epsilon} \exp(-H(s)/\epsilon)\,F(s)\,ds.
\end{align*}
\begin{lemma} \label{lem:limit exists} 
 Assume that $ \mathcal F :=(g_\epsilon)_{\epsilon>0}$ is a family of uniformly bounded solutions to \eqref{g-eq}. Then $\mathcal F$ is also equi-continuous in $\R \setminus \{ 0 \} $.
 In particular, up to a subsequence $(\epsilon_n)_n$, there exists $ \hat{ g}  \in C^0(\R)$  such that 
\begin{equation*}
    \lim_{ \epsilon_n \to 0 }  g_{\epsilon_n}  = \hat{g} \quad  \text{ locally uniformly in } x \in \R \setminus \{ 0 \}. 
\end{equation*}
 
\end{lemma}
\begin{proof}
Since  $\vert \vert g_\epsilon\vert\vert \le C$ we have that 
\begin{equation*}
    \vert\vert h_\epsilon(x)\vert\vert \le \vert\vert M(x)\vert\vert C \int_x^\infty (1/\epsilon) \exp((H(x)-H(y))/\epsilon)\,dy \leq \tilde{c}(x) , \quad x >0
\end{equation*}
where $\tilde{c}(x)$ is independent on $\epsilon$. As a consequence for every $x$ the function $h_\epsilon(x) $ is uniformly bounded as $\epsilon\to 0$. 
We can repeat the same type of argument for $x <0$. 
We deduce that, for any $\delta>0$, there is a constant $C(\delta)$ such that if $\vert x\vert\ge\delta$, then $\vert\vert h_\epsilon\vert\vert \le C(\delta)$. Hence $(g_\epsilon)_{\epsilon>0}$ is equi-Lipschitz and equi-continuous away from $0$.
\end{proof}

We now prove that the limit $\hat{g } $ of the sequence $(g_\epsilon)_{ \epsilon>0 }$ constructed in Lemma \ref{lem:limit exists} satisfies the limit equation \eqref{limit_eq}. 

\begin{lemma}\label{lem:limit_solves_limit_eq}
Assume that $\hat{g}$ is the function constructed in Lemma \ref{lem:limit exists}. Then $\hat{g} $ is the unique solution to \eqref{limit_eq}. 
\end{lemma}
\begin{proof} 
Since the derivative of the function  $g_\epsilon$ is now controlled, we  apply Laplace's  method (see \cite{bender1999advanced}) to evaluate the integral \eqref{eq:int for h_epsilon} defining $h_\epsilon$.  For $x\ge\delta>0$ we have that, up to the subsequence provided by Lemma \ref{lem:limit exists},
\begin{align*}
    h_\epsilon(x)&=-\int_x^\infty \frac{1}{\epsilon}\exp(H(x)-H(y))/\epsilon)M(y)g_\epsilon(y)\,dy\approx\frac{1}{H'(x)}M(x) \hat{g} (x), \text{ as } \epsilon \to 0.  
\end{align*}
More precisely, we have that for $x \geq \delta >0 $ 
\[ 
h_\epsilon (x) - \dfrac{1}{x } M(x) \hat{ g } (x) =: R_\epsilon(\delta ) \to 0 \text{ as } \epsilon \to 0
\] 
Repeating the argument for $x <0 $, we deduce that 
\begin{equation*}
    g_\epsilon'(x) =h_\epsilon (x) \to h (x)  \coloneqq \frac{1}{H'(x)}M(x)\hat{g}(x),\quad\text{locally uniformly on }\R\setminus\{0\}.
\end{equation*}
Since the $0\le \langle \mathbf{e}_k,g_\epsilon\rangle\le C$, the same holds for $\hat{g}$ and thus $\hat{g}$ is a bounded solution to $H'(x)g'=M(x)g$, with some initial datum $g(x_0)$ for $x_0\not=0$. By the uniqueness of bounded solutions to the limit equation, Proposition \ref{prop: existence and uniqueness solution limiting problem}, we must have $\hat{g}(0)\in\ker M(0)$ and it is uniquely determined by the conservation of mass.
Since $\hat{g}$ is unique, the above holds along any subsequence.
\end{proof} 
To upgrade this to globally uniform convergence, we have to consider also the boundary layer at $x=0$ where the above calculations do not apply.
Our approach will be to find an explicit formula for approximating the solution, at $0$ and $\pm\infty$, in terms of special functions. 
First, we consider the special case of $H(x)=x^2/2$ to introduce the main arguments. Afterwards, in Section \ref{subsec_general_H} we indicate the changes that are necessary in order to treat more general $H$. 
\begin{assumption}
Let $H(x)=x^2/2$.
\end{assumption}

We first assume that $M(x)\equiv M$ is diagonalisable and then extend the argument to matrices $M$ that have non-trivial Jordan blocks and afterwards to non-constant $M$. If $M$ is diagonalisable, we can change the basis on $\R^n$ to obtain $n$ equations that involve only one component of $g$
\begin{equation*}
    \epsilon g_i''- x g'_i +\lambda_i g_i=0. 
\end{equation*}
Writing $g_i(x)=\exp(x^2/(4\epsilon))\, q_i(x)$ and $z=x/\sqrt{\epsilon}$ we obtain Weber's differential equation (see \cite{NIST:DLMF}) for the components $q_i$
\begin{equation} \label{Weber}
    q''_{i}+\left((1/2+\lambda_i)- \frac{z^2}{4}\right)q=0.
\end{equation}
 This equation has well-known solutions in terms of the parabolic cylinder functions, cf., e.g., \parencite[Section 12]{NIST:DLMF} and \cite{pcf_paper}. It will be convenient to choose the two independent solutions by their behaviour at infinity. We write $b_i=1/2+\lambda_i$ and suppress the index $i$ in the following. The independent solutions that we use are
 \begin{equation*}
     \begin{aligned}
      U(z,b) &\sim \exp(-z^2/4)\, z^\lambda=\exp(-z^2/4)\, z^{b-1/2},\quad z\to \infty,\\
      V(z,b) &\sim \exp(z^2/4)\,z^{-\lambda-1}=\exp(z^2/4)\, z^{-b-1/2}\quad z\to\infty. 
     \end{aligned}
 \end{equation*}
 Recall that $\Re \lambda <0$ (whenever $\lambda\not=0$) and thus $\exp(z^2/4) U(z,b)$ decays algebraically.
 \begin{remark}\label{convention_pcf}
 In the literature for the Parabolic Cylinder Functions, the ODE \eqref{Weber} is usually written as 
 \begin{equation*}
     q''-(a+\frac{z^2}{4})q=0,\quad a\coloneqq-b
 \end{equation*}
 and accordingly $U(z,a)\sim \exp(-z^2/4)\,z^{-a-1/2}$ and $V(z,a)\sim \exp(z^2/4) \,z^{a-1/2}$. 
 \end{remark}
 \begin{lemma}\label{voc_pcf}
  Let $g$ be a bounded solution to 
  \begin{equation*}
      g''(z)- z g'(z) + \lambda g(z) + r(z)g(z)=0,
  \end{equation*}
  where $b=\lambda+1/2$, and $\Re (\lambda) <0$ and $\vert r(z) \vert \sim z^{-\alpha},\,\alpha>0$, as $z\to \infty$.\par
  Then, $\lim_{z\to\infty} g(z)=0$.
 \end{lemma}
 \begin{proof}
 Since $g$ is bounded, $r(z)g(z)=O(z^{-\alpha})$ as well and we write $R(z)\coloneqq r(z)g(z)$. We pick some sufficiently large starting point $z_0$ and write a variation of constants formula for $g$.
 For this we consider the Wronskian
 \begin{equation*}
     W(z)\coloneqq \exp(z^2/2) \det\begin{pmatrix} U &V\\U' &V' \end{pmatrix}.
 \end{equation*}
 As is easy to see from the equation, we have $W'(z)=z W(z)$, i.e., $W(z)=W_0\exp(z^2/2)$. By Cramer's rule we can write the general solution $g$ as
 \begin{align*}
     g &= \exp(z^2/4)\left( c_1 U(z,b) - U(z,b)\int_{z_0}^z \frac{\exp(s^2/4)}{W(s)} R(s) V(s,b)\,ds\right) \\
      &+ \exp(z^2/4)\left(c_2 V(z,b) + V(z,b)\int_{z_0}^z \frac{\exp(s^2/4)}{W(s)} R(s) U(s,b)\,ds\right).
 \end{align*}
 We estimate both summands of the right-hand side separately. For the first summand, we recall that $\exp(z^2/4)U(z,b)\sim z^{b-1/2}$, as $z\to \infty$. We have, after cancelling the exponentials
\begin{equation*}
     z^{b-1/2}\int_{z_0}^z \frac{1}{W(s)} R(s) V(s,b)\,ds \sim z^{b-1/2}\int_{z_0}^zs^{-b-1/2-\alpha}\,ds\sim z^{-\alpha},\quad z\to\infty.
\end{equation*}
For the exponentially growing solution $V$ we again have a compatibility condition coming from the boundedness of $g$. Since the other terms are already bounded, we must have
\begin{equation*}
 c_2+\int_{z_0}^\infty\frac{\exp(s^2/4)}{W(s)} R(s) U(s,b)\,ds =0, 
\end{equation*}
where the integral is obviously finite. Therefore we may write $g$ as
\begin{align*}
     g &= \exp(z^2/4)\left( c_1 U(z,b) - U(z,b)\int_{z_0}^z \frac{\exp(s^2/4)}{W(s)} R(s) V(s,b)\,ds\right) \\
     &- \exp(z^2/4)V(z,b)\int_{z}^\infty\frac{\exp(s^2/4)}{W(s)} R(s) U(s,b)\,ds.
\end{align*}
Using $V(z,b)U(z,b) \sim z^{-1}$ and $\exp(s^2/4)/W(s) = W_0^{-1} \exp(-z^2/4)$, as $z\to\infty$, we see that this last term also decays. 
 \end{proof}
 \begin{remark}
 If one repeats the proof in the original variable $x$, on obtains that $g\sim x^{-\alpha}$, as $x\to\infty$, independent of $\epsilon$.
 The above clearly remains valid if one replaces $r(z) g(z)$ by $a(z)g(z)+ b(z)$ where $a$ and $b$ both decay like a power law, and we still take $g$ bounded.
 In fact, the proof of this lemma even shows that in this case $g\sim z^{-\min(1,\alpha)}$, as $z\to\infty$.
 \end{remark}
 
\subsection{Behaviour of \texorpdfstring{$ g_\epsilon$, as $x \to \infty$}{∞}}
 Recall from Assumption \ref{assumption_M_structure}, that there is $\alpha>0$, such that $r(x)\coloneqq M(x)-M_\infty=O(x^{-\alpha})$, as $x\to\infty$. Thus, Lemma \ref{voc_pcf} applies and we can determine the behaviour of $g$ as $|x|\to\infty$.
 \begin{proposition}\label{pcf_g_infty}
 For $\epsilon>0 $ fixed, and $0\le g_\epsilon \le C$ bounded solution to \eqref{g-eq} with mass $m>0$, the limit 
 \begin{equation*}
     g_\infty \coloneqq\lim_{x\to\infty} g_\epsilon (x) 
 \end{equation*}
 exists and lies in $\ker M_+=\ker(\lim_{x\to\infty} M(x) )$. Further, $g_\infty$ is uniquely determined by $m$. The same is true at $x=-\infty$ with  $M(-\infty)$.
 \end{proposition}
 \begin{proof}
 We may assume that $M_\infty$ is in Jordan canonical form and the first block is the one-dimensional kernel. If we are in a Jordan block for some eigenvalue $\lambda\not=0$, then the equation for some component $g_i$ that belongs to a generalised eigenvector is only coupled to $g_{i+1}$.
 \begin{equation*}
     \epsilon g_i'' -x g_i' +\lambda  g_i +r g_i + g_{i+1}=0.
 \end{equation*}
 Write $z=x/\sqrt{\epsilon}$.
 If $k$ is the larges index belonging to the Jordan block, then we apply Lemma \ref{voc_pcf} to $g_k$ and see that $g_k(z)\to 0$. To be precise, the lemma gives that $g_k\sim z^{-\alpha}$ as $z\to\infty$, which is the same decay as $r(z)$ in the statement of Lemma \ref{voc_pcf}. Thus, we may apply this lemma inductively to $g_{k-i},\,i=1,\dots$ until we have exhausted the Jordan block.\par
 We conclude that all components except the first vanish as $x\to \infty$. Since the first Jordan block is spanned by the kernel, $g_1$ solves simply $\epsilon g_1''-x g_1'+ rg_1=0$. The two homogeneous solutions are now $C_1 \exp(x^2/(2 \epsilon))$ and the constant function $C_2$, the variation of constants argument is now much simpler than Lemma \ref{voc_pcf} (it works in a way very similar to Lemma \ref{prop: existence and uniqueness solution limiting problem}) and yields that in the limit $g_1$ behaves like a constant. Conservation of mass then determines the value of this constant. For $-\infty$, one changes coordinates to $y=-x$, which keeps the structure of the equation identical.
 \end{proof}
 \subsection{Effect of diffusion as \texorpdfstring{$x\to0$}{x→0}}\label{boundary_layer}
 Define $z=x/\sqrt{\epsilon}$, $X_\epsilon\coloneqq \epsilon^{1/4}$ and $Z_\epsilon \coloneqq \epsilon^{-1/4}$. Then $\vert x \vert \le X_\epsilon$ if and only if $\vert z\vert \le Z_\epsilon$.
 After bringing $M(0)$ to Jordan normal form and writing again $r(\bullet)=M(\bullet)-M(0)$ we have the same structure of equations as before
 \begin{equation}\label{eq:phi_pcf}
     \phi_i''-z \phi_i' +\lambda \phi_i + r \phi_i +\phi_{i+1}=0.
 \end{equation}
 Where now \begin{equation*}
     \phi_\epsilon(z)\coloneqq g_\epsilon(x/\sqrt{\epsilon})
 \end{equation*}
 \begin{lemma} 
 For every $\phi_\epsilon$ solutions of \eqref{eq:phi_pcf}, there exists a bounded solution   $\Psi_\epsilon=(\Psi_i)_{i=1}^n$ to 
  \begin{equation*}
      \Psi_\epsilon''-z \Psi_\epsilon' +M(0)\Psi=0,
  \end{equation*}
  such that for $\epsilon$ small enough
  \begin{equation*}
      \sup_{\vert z\vert \le Z_\epsilon/2^K}\vert \Psi_\epsilon(z)-\phi_\epsilon(z)\vert = O(\epsilon^{1/4}), \quad\text{as } \epsilon\to0.
  \end{equation*}
  Where $K$ is the size of the largest Jordan block of $M(0)$.
 \end{lemma}
 \begin{proof}
 A Taylor expansion around $0$ gives $M(x)=M(0)+M'(0)(\sqrt{\epsilon}x)+O(\epsilon)$. We see that for $\vert x\vert \le X_\epsilon$, we have $\vert R(z)\vert=\vert (M(z)-M(0))g\vert=O(\epsilon^{1/4})$. 
 Again, we bring $M(0)$ to Jordan normal form. As for the proof of Lemma \ref{voc_pcf}, consider the eigenvectors first and argue inductively for the rest of the Jordan blocks. Let $b=1/2+\lambda_i$, where $\lambda_i$ is an eigenvalue of $M(0)$.
 \begin{align*}
     \phi_i &= \exp(z^2/4)\left( c_1(\epsilon) U(z,b) - U(z,b)\int_{z_0}^z \frac{\exp(s^2/4)}{W(s)} R(s) V(s,b)\,ds\right) \\
      &+ \exp(z^2/4)V(z,b)\left( c_2(\epsilon) + \int_{z_0}^z \frac{\exp(s^2/4)}{W(s)} R(s) U(s,b)\,ds\right).
 \end{align*}
 In the same way as in the proof of Lemma \ref{voc_pcf}, the first summand is seen to be $\sim c_1(\epsilon) z^\lambda+ O(\epsilon^{1/4})$, as $\epsilon\to0$, for $z\in [-Z_\epsilon,Z_\epsilon]$. Since the $\phi_\epsilon$ are uniformly bounded, there must be a $C>0$ such that 
 \begin{equation*}
    \vert \kappa_\epsilon \vert \coloneqq \vert c_2 +\int_0^{Z_\epsilon} \frac{\exp(s^2/4)}{W(s)}R(s) U(s,b)\,ds\vert\le C\left(\exp(Z_\epsilon^2/4)V(Z_\epsilon,b)\right)^{-1}.
 \end{equation*}
 (The integrand satisfies $\vert\frac{\exp(s^2/4)}{W(s)}R(s) U(s,b)\vert\le  C\epsilon^{1/4}\exp(-s^2/2)s^\lambda$ and $c_2(\epsilon) =o(\epsilon^{1/4})$, as $\epsilon\to0$, follows.)
 We write now
 \begin{align*}
     \phi_i &= \exp(z^2/4)\left( c_1(\epsilon) U(z,b) - U(z,b)\int_{z_0}^z \frac{\exp(s^2/4)}{W(s)} R(s) V(s,b)\,ds\right) \\
      &+ \exp(z^2/4)V(z,b)\left( \kappa_\epsilon - \int_{z}^{Z_\epsilon} \frac{\exp(s^2/4)}{W(s)} R(s) U(s,b)\,ds\right)\\
      &\coloneqq \Psi_\epsilon^{(i)}(z) +  \kappa_\epsilon V(z,b)\exp(z^2/4) + U(z,b)\int_{z_0}^z \frac{\exp(s^2/4)}{W(s)} R(s) V(s,b)\,ds
 \end{align*}
 Here, $\Psi$ is the solution we seek: For $\epsilon$ small enough, the first remainder term will be exponentially small for, e.g., $\vert z\vert \le Z_\epsilon/2$ while the second is $O(\epsilon^{1/4})$, as $\epsilon\to0$, as $\epsilon\to0$. We conclude by repeating this inductively along the Jordan blocks. Notice that in each step we half the size of the interval again, leading to $\vert z\vert\le Z_\epsilon/2^K$ in the end.
 \end{proof}
 \begin{lemma}
 Let $\phi_\epsilon$ be a solution to equation \eqref{eq:phi_pcf}, as in the Lemma above. Set $a_\epsilon(z)\coloneqq \langle \mathbf{1},\phi(z)\rangle$ and $m_\epsilon=\sum_{k=1}^n \Psi_k$, then $a_\epsilon(z)\equiv a_\epsilon(0)=m_\epsilon+o(\epsilon^{1/4})$, as $\epsilon\to0$. 
 \end{lemma}
 \begin{proof}
The function $a_\epsilon$ solves $0=a''(z)-z a'(z)$ since $\mathbf{1}\in \ker M(z)$ for all $z$. Thus, $a_\epsilon$ is either  unbounded or constant. We conclude it is constant by the mass constraint. Next, we use that $\phi$ also has conservation of mass and on the interval $[-Z_\epsilon/2^K,Z_\epsilon/2^K]$ the Lemma above implies the second equality. 
 \end{proof}
 Recall the following from above.
 \begin{lemma}
  For $m>0$, let $g(x)$ be the unique bounded solution to the first-order limit equation \eqref{limit_eq}, then $g(0)\in \ker M(0)$.
 \end{lemma}
 \begin{lemma}\label{pcf_g_0}
  For $\vert z\vert \le Z_\epsilon$, respective $\vert x\vert \le X_\epsilon$, 
  \begin{equation*}
      \vert\vert\phi_\epsilon(z)-g(0) \vert\vert=\vert\vert g_\epsilon(x)-g(0) \vert\vert =O(\epsilon^{1/4}), \quad \epsilon\to0.
  \end{equation*}
 \end{lemma}
 \begin{proof}
For $M(\bullet)\equiv M(0)$, we can immediately see the uniqueness of the bounded solutions $\Psi$. Indeed, we put $M(0)$ again into the (complex) Jordan canonical form, and argue inductively over over the Jordan blocks. The claim we need to show is, that for $\Re \lambda<0$, the only bounded (on $\R$)  solution $q$ of
\begin{equation*}
     q''+\left((1/2+\lambda)- \frac{z^2}{4}\right)q=0,
\end{equation*}
that can give a bounded $\Psi_i=\exp(z^2/4)q(z)$, is $q\equiv 0$. Since $q$ will be a combination of the solutions $V(z, -1/2-\lambda),\, U(z,-1/2-\lambda)$, this claim follows from the connection formulas:
\begin{align}
    U(-z,a)&=-\sin(\pi (a)U(z,a)+\frac{\pi}{\Gamma(\frac{1}{2}+(a))}V(z, a),\\
    V(-z, a)&=\sin(\pi (a))V(z, a) +\frac{\cos(\pi a)}{\Gamma(1/2-(a))} U(z, a),
\end{align}
with $a\coloneqq -1/2-\lambda$, which implies that the decaying solutions as $z\to+\infty$, pick up exponential growth at $z\to-\infty$ and vice versa. Lastly, the kernel-component of $\Psi$ is constant and determined by the mass.
By definition, $\Psi_\epsilon$ is bounded on the whole of $\R$.
 On the other hand, letting $\ell \in \ker M(0)$, such that $\sum_k \ell_k=1$, we have another solution given by $m_\epsilon\,\ell$ of the same equation. By uniqueness, this solution must equal $\Psi_\epsilon$. Also, we have $m_\epsilon \to m$, since $m:=\sum\phi_k=\sum_k \Psi_k+O(\epsilon^{1/4})$, as $\epsilon\to0$, for $\vert z\vert \le Z_\epsilon/2^K$.
 \end{proof}
We combine the results of this section to obtain the uniform convergence on $\R$, i.e. to prove Theorem \ref{thm: convergence as epsilon to zero}. 
 \begin{proof}[Proof of Theorem \ref{thm: convergence as epsilon to zero}]
 This is obtained by combining the local uniform convergence result from Lemma \ref{lem:limit_solves_limit_eq} with Lemma \ref{pcf_g_infty} and Lemma \ref{pcf_g_0}. The only addition needed is the observation that the latter result holds for $\vert x\vert\le C \epsilon^{1/4}\gg\sqrt{\epsilon}$, which is the scale where the Laplace approximation of $g'_\epsilon(x)$ breaks down.
 \end{proof}

 \subsection{General \texorpdfstring{$H$}{H}}\label{subsec_general_H}
  Lastly, we indicate the necessary minor modifications to the above calculations that allow to consider more general $H$.
 For now, we focus on the boundary layer at $0$, to treat the behaviour at $\infty$ with the same tools, we have to make additional assumptions.
 We still consider $g=\exp(H(x)/(2\epsilon))q(x)$ and write $x=z\sqrt{\epsilon}$. 
 \begin{equation}\label{pcf_general_H}
    \begin{aligned}
     0&=q''(z)+( \frac{H''(0)}{2}-\lambda_i-\frac{H'(0)^2\,z^2}{4})q(z)+ R(z)+M(z)-M(0),\\
     R(z)&=(H''(z)-H''(0))/2-[(H'(z))^2-H'(0)^2 z^2]/4
 \end{aligned}
 \end{equation}
 Since $H''(0)>0$ the structure of the equation is as before and the remainder $R(z)+M(z)-M(0)$ still is $O(\epsilon^{1/4})$ for $\vert z\vert \le Z_\epsilon\coloneqq \epsilon^{-1/4}$ as before.\par
 At $z=\infty$, however, we need to explicitly make the assumption that $R(z)$ is bounded like this. In particular we have to assume that the limits of $H'/z$ and $H''$ exist, to even make sense of the approximate form \eqref{pcf_general_H}. 
 In this way, we obtain the uniform convergence of Theorem \ref{thm: convergence as epsilon to zero} also for general $H$.
 
 \section{Constant flux solutions}\label{sec_fluxes}
 Since the system that we studied does not satisfy the detailed balance property, the steady states are non-equilibrium steady states and we expect to have non zero fluxes at the steady state. In order to see this, we explain in this section how one can define the fluxes for system \eqref{disc_cont_eq}. 
Recall Assumption \ref{assumption_H} and that $\ell(x)\in \ker M(x)$ is the unique kernel element such that $\mathbf{1} \cdot \ell(x)=1$.
Notice that if we assume $M(x) \equiv 0 $ then the definition of the fluxes is straightforward, since the solution will be of the form $F(t,x)\mathbf{1}$ and only depend on $x$. Every scalar Fokker-Planck equation can be written as a continuity equation 
\begin{equation}\label{flux_sign-convention}
    \partial_t-\partial_x(\epsilon \partial_x f +H'f)=\partial_t f -\Div(J_x)=0,
\end{equation}
where $J_x=J_x(f)$ is the flux in the $x$-variable. Accordingly, we define
 \begin{equation*}
    J(f) \coloneqq \epsilon -\partial_x f - H'\,f=-\epsilon\partial_x g\,\gamma.
 \end{equation*}
 The detailed balance property for one-dimensional Fokker-Planck equations (without periodic boundary conditions) is always satisfied and indeed the stationary solutions have always zero fluxes. 
 While this is always the case for scalar problems, it is also easy to see that in the present case, $J(f)=J(g \gamma)$ only vanishes identically if $g$ is constant. In particular, this cannot happen if $M(x)$ is not independent of $x$. 

For $\epsilon\ll1$ however, we can at least show that detailed balance holds approximately.
For vanishing $\epsilon$, the stationary solution $g$ converges to the unique bounded solution to $\bar{g}'=(1/H') M \bar{g}$, which has the property that $\bar{g}(0)\in\ker M(0)$ and thus $\bar{g}'(0)=0$. 
In Section \ref{sec_asymptotics}, we have shown that this implies $\bar{g}'(0)=0$, i.e.,
\begin{equation*}
     J(\bar{g}\gamma)\sim\sqrt{\epsilon}\,\delta_0,\quad \epsilon\to0,
\end{equation*}
where $\delta_0$ denotes the Dirac delta in $x=0$.
For calculating $\epsilon\frac{d}{dx}g_\epsilon$, we use the integral representation of $h_\epsilon=g_\epsilon'$ from the proof of Lemma \ref{lem:limit_solves_limit_eq} and the uniform convergence $g_\epsilon\to\bar{g}$ of said lemma. We have
\begin{align*}
   \epsilon \frac{d}{dx}g_\epsilon(x) &=-\epsilon\int_x^\infty \frac{1}{\epsilon}\exp(H(x)-H(y))/\epsilon)M(y)g_\epsilon(y)\,dy\\
   &\approx\epsilon\frac{1}{H'(x)}M(x)\bar{g} (x)=o_\epsilon(1),\quad \epsilon \to 0. 
\end{align*}
Hence, we also obtain
 \begin{equation*}
     J(g_\epsilon\gamma)=o_\epsilon(1)\,\delta_0, \quad \epsilon \to 0.
 \end{equation*}
 Recall that the chemical system made by the macroscopic chemicals $A$ and $B$ and by the chemical reactions $\rho_0$ and $- \rho_0$ satisfies the detailed balance property. However, this macroscopic chemical system interacts with $K $. This keeps the dynamics of the macroscopic chemicals (slightly) out of equilibrium. 
 
The analogous ``discrete flux'' is harder to define and requires a discrete analogue of the divergence. In the continuous variable we consider a scalar quantity and the flux is a vector field, which describes the direction in which the transport of said quantity takes places. In analogy, in the discrete variable we consider a vector, hence the flux $I$ will be a matrix, which we should accordingly read as a function $I:\{1,\dots,n\} \to \R^n$ with $I_{i,j}=I(i)_j$.
Consider a general linear reaction network of $n$ species that conserves the total mass. Let $M=(m_{i,j})_{i,j=1}^n$ be the matrix of reaction rates $(i)\to (j)$. If at some point the concentrations are in state $v\in \R^n$, we consider the net conversion of species $(i)$ into $(j)$ in that instant. 
\begin{definition}
Define the linear map $F_M:\R^n\to \R^{n\times n}$, via
\begin{equation*}
    (F^M\,v)_{i,j}\coloneqq \begin{cases}
       m_{i,j}v_i-m_{j,i}v_j,\quad &i\not=j\\
       0,\quad &i=j.
    \end{cases} 
\end{equation*}
\end{definition}
We recover the usual species formation function starting from $F_M$ via 
\begin{equation*}
M\,v =(F^M\,v)\,\mathbf{1}.
\end{equation*}
By construction we have $F_M(v)_{i,j}=-F_M(v)_{j,i}$, and clearly $F_M(v)$ being skew-symmetric is just conservation of mass. 
For now it will be more convenient to write the stationary solutions of \eqref{disc_cont_eq} as $f(x,k)$ for $(x,k)\in \R\times\{1,\dots, n\}$. To define the flux at a point $(x,k)$, we need to choose an orientation on the reactions in the network. Here we simply take the reaction $k\to\ell$ as positively oriented if $k<\ell$ and $\ell\to k$ to be negatively oriented. In view of the sign convention of fluxes in the continuous variable, i.e., \eqref{flux_sign-convention}, we make the following definition of the flux 
\begin{equation}\label{eq:flux_def}
    \mathcal{J}(\bar f)(x,k) \coloneqq (J_\epsilon(\bar{f}_\epsilon),\mathbf{0})^\dagger - \sum_{\ell\not= k} \sgn(\ell-k)F_M(\bar{f}_\epsilon)(x)_{k,\ell}(0,\mathbf{e}_\ell)^\dagger.
\end{equation}
where $\mathbf{0}, \mathbf{e}_k$ denote the zero vector and the standard basis vectors, respectively.
Since detailed balance means that there exists $v\in \R^n$ such that $m_{i,j} v_i=m_{j,i}v_j$ and such $v$ lies in $\ker M$, we can directly relate detailed balance to $I_M$.
 \begin{lemma}
 For a connected network 
 \begin{equation*}
     \begin{tikzcd}
	 {(i)} & {(j)}
	\arrow["{\rho_{i,j}}", shift left=2, from=1-1, to=1-2]
	\arrow["{-\rho_{j,i}}", shift left=2, from=1-2, to=1-1]
\end{tikzcd}
 \end{equation*}
 detailed balance holds if and only if $F_M(v)=0$ or equivalently $I_M(v)^\dagger =I_M(v)$ hold for $v\in \ker M$.
 \end{lemma}
Therefore, the structure of the fluxes is precisely as described in Theorem \ref{thm:fluxes_main results}. 

\begin{remark}
For a vector field $b$, we can write the divergence as $\Div b=\tr \nabla b=\sum_k \partial_k b$. We try to reproduce this structure. In the discrete setting the analogue of differentiation should be a difference, hence we expect to have
\begin{equation*}
    Mv= F_m(v)\,\mathbf{1} = \delta(I_M(v))\,\mathbf{1},
\end{equation*}
where $I$ is now the discrete flux, $\delta(\bullet)$ is a finite difference operator and multiplication by $\mathbf{1}$ provides the summation over the finite differences.
We define 
\begin{equation*}
    \delta: \R^{n\times n}\to \R^{n\times n}: (a_{ij})_{i,j=1}^n\mapsto (\delta(a)_{i,j=1})^n \coloneqq (a_{i,j}-a_{j,i})_{i,j=1}^n
\end{equation*}
and 
\begin{equation*}
    I_M(v)\coloneqq (m_{i,j}v_i)_{i,j}^n
\end{equation*}
 For $\epsilon\to 0$ we have
\begin{equation*}
    I_{M(x)} g\,\gamma \approx m I_{M(0)}(\ell(0)) \gamma
\end{equation*}

If we specialise again to the $M(x)$ as in \eqref{intro_disc_cont_eq}, and if we additionally assume $a_k(x)=b_k(x)$, and thus detailed balance in the reactions \eqref{network}, then we can write $Mv$ as a second order (central) finite difference.
We define the forward and backward differences as 
\begin{align*}
    \delta^+ (v)&\coloneqq\sum_k (v_{k+1}-v_k)\,\mathbf{e}_k\\
    \delta^- (v) &\coloneqq\sum_k (v_k-v_{k-1})\,\mathbf{e}_k
\end{align*}
and accordingly we have that 
\begin{equation*}
    M f(t,x)= \delta^-(\delta^+( a f )),
\end{equation*}
with the point-wise product $av=\sum_k (a_k f_k)\, \mathbf{e}_k$. Informally, the term $M\,f$ is a discretisation of a second-order derivative of the form $\partial_k^2 (a(x,k)\,f(t,x,k))$. This explains why solutions to the limit equation $g'(y)=M(y)g(y)$ in Section \ref{sec_asymptotics} behave quite similarly to a parabolic PDE, with time variable $y$ and spatial variable $k$. On the other hand, for general $M$, and general networks, it is not possible to write $Mf$ as a composition of a forward and a backward difference operator.
 \end{remark}
 
 \bigskip 
 
 \textbf{Acknowledgements.} The authors gratefully acknowledge the financial support of the Deutsche
Forschungsgemeinschaft (DFG, German Research Foundation) through the collaborative re-
search centre ”Analysis of criticality: from complex phenomena to models and estimates” (CRC
1720, Project-ID 539309657) and the Bonn International Graduate School of Mathematics at the
Hausdorff Center for Mathematics (EXC 2047/2 Project-ID 390685813).
The funders had no role in study design, analysis, decision to publish, or preparation of the manuscript.

\bigskip

\textbf{Data Availability.} No data are associated with this article.
\bigskip

\textbf{Statements and Declarations.}
The authors declare no competing interests.  

 \printbibliography
\end{document}